\documentclass[12pt,a4paper]{article}
\usepackage{amssymb,amsbsy,amsmath,amsfonts,amssymb,amscd}
\usepackage{latexsym}

\input xy
\xyoption{all}

\newtheorem{theo}{Theorem}[section]
\newtheorem{remarkk}[theo]{Remark}
\newenvironment{rem}{\begin{remarkk}\rm}{\end{remarkk}}
\newtheorem{definition}[theo]{Definition}
\newenvironment{devery easyfi}{\begin{definition}\rm}{\end{definition}}
\newtheorem{prop}[theo] {Proposition}

\newtheorem{lemma}[theo]{Lemma}
\newtheorem{example}[theo]{Example}

\newtheorem{comp}{Computational Fact}[section]

\newcommand{\R}{\ensuremath{\mathbb{R}}}

\newcommand{\NN}{\ensuremath{\mathbb{N}}}

\newcommand{\HH}{\ensuremath{\mathbb{H}}}
\newcommand{\PP}{\ensuremath{\mathbb{P}}}

\def\Bbb{\bf}

\def\C{{\Bbb C}}

\def\Z{{\rm Z}}

\def\SmG{{\rm SmallGroup}}
\def \Mid{\quad \vrule \quad}

\def\BP{{\hbox{\ensuremath{\mathbb{P}}}}}

\def\BN{{\hbox{\ensuremath{\mathbb{N}}}}}
\def\BF{{\hbox{\ensuremath{\mathbb{F}}}}}
\def\BT{{\hbox{\ensuremath{\mathbb{T}}}}}
\def\PSL{{\hbox{\rm PSL}}}
\font \msbm= msbm10
\def\rfish{\mathbin {\hbox {\msbm \char'157}}}

{ \noindent {\em Proof: }}%
{\hspace*{\fill}$\Box$\par \medskip }
\newcommand{\Proof}{{\it Proof. }}
\newcommand{\qed}{\hspace*{\fill}$Q.E.D.$}

\begin{document}
\title{The classification of surfaces  with $p_g = q = 0$ isogenous to a product of curves}
\author{I. C. Bauer, F. Catanese, F. Grunewald \\
\small Mathematisches Institut der Universit\"at Bayreuth \\
\small Universit\"atsstr. 30 \\
\small 95447 Bayreuth \\
\small Mathematisches Institut der Heinrich-Heine-Universit\"at D\"usseldorf\\
\small Universit\"atsstr.  \\
\small 40225 D\"usseldorf }
\date{\today}
\maketitle

\tableofcontents

\newpage
%%%%%%%%%%%%%%%%%%%%%%%%%%%%%%%%%%%%%%%
\section{Introduction}
%%%%%%%%%%%%%%%%%%%%%%%%%%%%%%%%%%%%%%%

It is well known that an algebraic curve of genus zero is isomorphic 
to the projective line. The
search for an analogous statement in the case of algebraic  surfaces 
led Max Noether to
conjecture that a smooth regular (i.e.,
$q(S) = 0$) algebraic surface with vanishing geometric genus ($p_g(S) 
= 0$) should be a rational
surface. The first counterexample to this conjecture was provided by 
Federigo Enriques in 1896 (\cite{enr96}, and also \cite{enrMS},
I, page 294), who introduced the so called Enriques surfaces by considering the 
normalization of sextic
surfaces in 3-space double along the edges of a tetrahedron.
Enriques surfaces are of special type, 
   nowadays a large number of surfaces of general type with $p_g = q = 
0$ is known, but the first
ones  were constructed in the thirties  by   Luigi Campedelli and Lucien 
Godeaux (cf. \cite{Cam},
\cite{god}: in their honour  minimal surfaces  of general type with 
$K^2 = 1$  are called
numerical Godeaux surfaces, and those with
$K^2 = 2$ are called numerical  Campedelli surfaces).

In the seventies, after rediscoveries of these old examples, many new 
ones were found through
the efforts of several authors (cf. \cite{bpv}, pages 234-237 and 
references therein). In
particular, in the spirit of Godeaux' method to produce interesting 
surfaces as quotients $S =
Z/G$ of simpler surfaces by the free action of a finite group $G$, 
Arnaud Beauville proposed a very
simple construction  by taking as $Z$ the product $Z=C_1 \times 
C_2$ of two curves of respective genera  
$g_1 := g(C_1), g_2 : = g(C_2) \geq 2$, 
together with an action of a group $G$ of order 
$(g_1 - 1 )(g_2 -1 )$ (this method
produces surfaces with
$K^2 = 8$). 
He also gave an explicit example as quotient of two Fermat curves
(in \cite{BaCa} it was shown that his example leads to exactly two non 
isomorphic surfaces).

Generalising Beauville's construction we study here surfaces 
 $S$ isogenous to a
product of two curves, i.e., surfaces which have a finite unramified cover
which is biholomorphic to a product of two curves.  One says that 
the surface $S$ is isogenous to a higher product if both curves have genus bigger or equal
to $2$ (this condition is equivalent to $S$ being of general type). 

It turns out that any surface with $p_g = q = 
0$ and isogenous to a  product is either $\PP^1 \times \PP^1$ or it is isogenous to a higher product
(this happens since $\chi (S) : = \chi(\mathcal{O}_S) =1 \Longrightarrow \chi (C_1 \times C_2) = (g_1 -1) (g_2 -1) >
0$, whence either both $g_i$'s are $\geq 2$, or both  $g_i$'s are $ = 0$).

By results of \cite{cat00} any
 surface $S$ isogenous to a higher product has a 
unique minimal realisation $S\cong (C_1\times C_2) /G$ where
$G$ is a finite group acting freely on $C_1\times C_2$ and with $g_1 := g(C_1), g_2 : = g(C_2) \geq 2$
chosen as small as possible. 
The action of $G$ can be seen to respect the product structure of $C_1\times C_2$. This
means that there  are the following
two possibilities. Either there are actions of $G$ on $C_1$ and $C_2$ such
that the action of $G$ on $C_1\times C_2$ is the diagonal action,and  if this happens we
speak of the {\em unmixed} case. 
Or there are elements in $G$ which interchange 
$C_1$ and $C_2$, and if this happens we
speak of the {\em mixed} case. Obviously, in the mixed case 
$C_1$ and $C_2$ have to be
biholomorphic to each other.   
 
In this paper we carry out the classification of all 
smooth projective surfaces $S$ isogenous to a product
with $p_g(S)=q (S) = 0$. Note that if $S$ is of general type $p_g = 0$ implies  $q = 0$, since for a
surface of general type $\chi (S) : = \chi(\mathcal{O}_S) = 1 + p_g -q \geq 1$. 

We can henceforth assume without loss of generality that $ S \ncong \PP^1 \times \PP^1$ 
and therefore that $S$ is of general type.

First invariants of such surfaces are the group $G$ of the minimal
realisation $S\cong (C_1\times C_2) /G$ and the genera of $C_1$ and $C_2$.

It turns out that the surfaces $S$ which can be obtained 
for a fixed finite group $G$ and with fixed genera $g_1 := g(C_1), g_2 : = g(C_2)$
fill out a finite number $N$ of irreducible connected components in  
the moduli space $\mathfrak{M}_{(1,8)}$ of 
minimal smooth complex
projective surfaces with  $\chi (S) = 1$  and $K_S^ 2 = 8$.
These turn out a posteriori to have the same dimension $D$. 

Our main result is:

\begin{theo}\label{theomai}
If $S$ is a smooth projective surface isogenous to a product
with $p_g(S)=q(S)=0$ 
and
with minimal realisation $S\cong (C_1\times C_2) /G$ then 
either $G$ is trivial and $ S \cong \PP^1 \times \PP^1$ or $G$ is one of the
groups in the following table and the genera of the curves $C_1,\, C_2$ are as listed   
in the table. The numbers of components $N$ in $\mathfrak{M}_{(1,8)}$ and
their dimension is given in the remaining two columns. 
\begin{center}
\begin{tabular}{|c|c|c|c|c|c|c|}
\hline
$G$  & $|G|$ & {\rm Type} & $g(C_1)$ & $g(C_2)$ & $N$ & $D$ \\
\hline
${\mathfrak A}_5$ & $60$ & {\rm unmixed} & $20$ & $3$  & $1$ & $1$ \\
${\mathfrak A}_5$ & $60$ & {\rm unmixed} & $5$  & $12$ & $1$ & $1$ \\
${\mathfrak A}_5$ & $60$ & {\rm unmixed} & $15$ & $4$ & $1$ & $1$ \\
${\mathfrak S}_4\times \Z_2$ & $48$ & {\rm unmixed} & $24$ & $2$  & $1$ & $3$
\\
${\rm G}(32)$ & $32$ & {\rm unmixed} & $4$ & $8$ & $1$ & $2$ \\
$\Z_5^2$ & $25$ & {\rm unmixed} & $5$ & $5$ & $2$ & $0$ \\
${\mathfrak S}_4$ & $24$ & {\rm unmixed} & $12$ & $2$ & $1$ & $3$ \\
${\rm G}(16)$ & $16$ & {\rm unmixed} & $4$ & $4$ & $1$ & $2$ \\
${\rm D}_4\times\Z_2$ & $16$ & {\rm unmixed} & $8$ & $2$ & $1$ & $4$ \\
$\Z_2^4$ & $16$ & {\rm unmixed} & $4$ & $4$ & $1$ & $4$ \\
$\Z_3^2$ & $9$ & {\rm unmixed} & $3$ & $3$ & $1$ & $2$ \\ 
$\Z_2^3$ & $8$ & {\rm unmixed} & $4$ & $2$ & $1$ & $5$ \\ 
${\rm G}(256,1)$ & $256$ & {\rm mixed} & $16$ & $16$ & $3$ & $0$ \\ 
${\rm G}(256,2)$ & $256$ & {\rm mixed} & $16$ & $16$ & $1$ & $0$ \\ 
\hline
\end{tabular}
\end{center}
Here ${\mathfrak A}_5$ is the alternating group on $5$ letters, $\mathfrak S_4$ is the
symmetric group on $4$ letters, ${\rm D}_4$ is the dihedral group of order $8$,
$\Z_n$ is the cyclic group of order $n$, ${\rm G}(32)$ and ${\rm G}(16)$ are
two groups of respective orders $32$, $16$ described in Sections \ref{suse44} and \ref{concrete}
and ${\rm G}(256,1)$,
${\rm G}(256,2)$ are two groups of order $256$ described in Sections
\ref{suse256} and \ref{concrete}.  
\end{theo}
\medskip

We see our main result as the solution in a very special case 
to the open problem that David 
Mumford set forth at the
Montreal Conference in 1980: "Can a computer classify all surfaces 
of general type with $p_g=0$? 
Our purpose is to show how computationally complex this question is and that
probably computers are needed even
if one asks a more restricted question. 

All known surfaces of general type  with $p_g = 0,\  
  K^2 = 8$ are quotients
$\HH \times \HH / \Gamma$ of the product of two upper half planes by 
a discrete cocompact subgroup of $\PSL(2,\R)\times \PSL(2,\R)$.
There are also quotients which are not 
related to products of curves; constructions of such surfaces using
quaternion algebras have been known since long, see for example
\cite{kug}, \cite{shav}. Also in this case a complete classification is 
possible. We shall elaborate on this in a forthcoming paper.

It is an interesting question whether there do exist surfaces of general type  with $p_g = 0, 
K^2 = 8$ which are not  quotients
$\HH \times \HH / \Gamma$ as above (observe that for $p_g = 0, 
K^2 = 9$ the universal cover is the unit disk in $\C^2$ by Yau's theorem
\cite{yau}). In particular, it is attributed to Hirzebruch the question of
existence of such surfaces of general type which are simply connected (they 
would be homeomorphic
to $\PP^1 \times \PP^1$ but not diffeomorphic, see \cite{free} and
\cite{don}).

Surfaces with
$p_g = q = 0$ were  also investigated  from other points of view. We 
would like to mention
several articles by M. Mendes Lopes and R. Pardini ( \cite{PardDP},
\cite{MLP1}, \cite{MLP2}) where the authors encountered them
in the course of studying the
failure of  birationality of the bicanonical map.

Concrete examples of rigid surfaces isogenous to a product 
have been given in
\cite{BCG}. 

Previously, in \cite{BaCa} the first two authors classified all 
smooth algebraic surfaces isogenous to a product and of unmixed type
with $p_g = q = 0$, and with $G$  a finite abelian group. 

They also gave a complete description of the connected components 
of the moduli space that
arise from these surfaces.

In this article we complete this classification admitting arbitrary
groups and treating also the mixed type.

While describing the organisation of the paper we shall now explain the steps
of our classification procedure in more detail. 

A family of surfaces isogenous to a product with associated group $G$ 
and with $ q=0$ is determined in the
unmixed case by a set of data which we
call a {\em ramification structure}. It consists of a pair of spherical systems of
generators 
$[g_{(1,1)}, \ldots, g_{(1,r)}]$, $[g_{(2,1)}, \ldots, g_{(2,s)}]$
for the group $G$ (i.e., a system of generators whose product equals the
identity), which are 'disjoint' in the sense that the union of the conjugates of
the cyclic subgroups generated by $g_{(1,1)}, \ldots, g_{(1,r)}$, resp.
$g_{(2,1)}, \ldots, g_{(2,s)}$ have trivial intersection. 
We exploit also the fact that
the geometric conditions impose very strong restrictions on the possible orders
of the elements $g_{(i,j)}$.

 We are able to classify these by combinatorial
methods of finite group theory. 
 Riemann's  existence theorem guarantees in fact that
for any  ramification structure there
is an irreducible family of surfaces isogenous to a product 
with the given ramification structure.  

In the mixed case we follow a similar approach.

In Section \ref{sec1} we fix the algebraic set up and classify all the possible
types (i.e., tuples of orders) of the spherical systems of generators. In fact,
the conditions on the possible tuples of orders are strong enough to leave
only finitely many possibilities also for the orders of the finite groups
which have to be considered.

In Section \ref{groups}
we introduce the action of the product of the braid group with
${\rm Aut}(G)$ on the set of 'disjoint' 
spherical systems of generators, which  reflects the
deformation equivalence of the associated surfaces.

In  Section \ref{basi} we collect some basic 
results on surfaces isogenous to a
product and show how they correspond to the  algebraic
data introduced previously.

In Sections \ref{classiumi}, \ref{classimi} we carry out 
the complete classification of
all finite groups occuring as groups $G$ associated to a surface isogenous to
a product with $p_g = 0$.
The procedure is simple: using the libraries of the MAGMA computer algebra
system \cite{MA}, which include all groups of 
order less then $2000$ (with the exception
of $1024$), we try to inspect all groups whose orders appear in the list 
obtained in Section
\ref{sec1} asking for the existence of suitable systems of generators. 
This turned
out not to be an easy task for two reasons: first, the number of groups which has to
be checked is much too high to be feasible to a direct computer calculation, 
second, the orders of the groups in question may be too high to be contained in the
standard group libraries. In order to prove our main result we have then to use
direct arguments (exploiting e.g. the solvability of groups whose order admits
only two prime factors), which allow us to reduce the
cardinality of the finite groups under consideration until we reach a region which
is covered by 
the MAGMA-library of small groups. 

We have not tried to minimise the amount of computer calculations needed. 
But we have tried to keep the complexity and time requirements for each single
calculation as small as possible.  Sections  
\ref{classiumi}, \ref{classimi} are to be understood 
as a Leitfaden through a maze of
little facts about finite groups.

After finishing the calculations we realized that (with some effort) all
computer calculations could be eliminated to give a, in fact much longer,
``hand made'' proof of our main results. 
We believe that the interest of our paper is twofold:
first of all the list of surfaces in Theorem \ref{theomai} contains many new
and interesting examples. Finding them is difficult, but establishing their
existence is easy. In particular, we devote Section
\ref{concrete} to a simple description of the groups and 
ramification structures occurring.
We hope that this description may be useful for working explicitly
with our surfaces. Secondly, it seems
interesting to us that it is at all possible to carry out a subcase of the Mumford
classification program. The fact that this is only a subcase
was the reason for us not to analyse our results further and free them
from computer calculations.

Finally, in section  \ref{modu} we calculate the number of orbits of the direct
product of the braid group with ${\rm Aut}(G)$ acting 
on the set of disjoint pairs
of spherical systems of generators. By this procedure we determine the exact
structure of the corresponding subset of the moduli space 
corresponding to  surfaces isogenous to a
product with
$p_g = 0$: in particular we determine the number  of irreducible
connected components and their respective  dimensions.

%%%%%%%%%%%%%%%%%%%%%%%%%%%%%%%%%%%%%%%%%%%%%%%%%%%%%%%%%
\section{Combinatorial preliminaries}\label{sec1}
%%%%%%%%%%%%%%%%%%%%%%%%%%%%%%%%%%%%%%%%%%%%%%%%%%%%%%%%%

This section contains simple combinatorial results
which are important as a first step  in the solution of 
the algebraic problem to which our classification can be reduced.
We also fix certain terminologies to be used later. 
The reader who finds   these preliminaries too dry to  swallow
might first want to read the subsequent section \ref{basi},
explaining how we pass from geometry to algebra.

%%%%%%%%%%%%%%%%%%%%%%%%%%%%%%%%%%%%%%%%%%%%%%%%%%%%%%%%%
\subsection{Group theoretic terminology}\label{groups}
%%%%%%%%%%%%%%%%%%%%%%%%%%%%%%%%%%%%%%%%%%%%%%%%%%%%%%%%%

Let $G$ be a group and $r\in\NN$ with $r\ge 2$. An $r$-tuple 
$T=[g_1,\ldots,g_r]$ of elements of $G$ is called a 
{\it spherical system of generators of $G$} if 
$g_1,\ldots,g_r$ is a system of generators of $G$
(i.e., $G=\langle\, g_1,\ldots,g_r\, \rangle$)
and we additionally have $g_1\cdot\ldots\cdot g_r=1$.

We call $r=:\ell(T)$ the length of $T$.

If $T=[g_1,\ldots,g_r]$ is an $r$-tuple of 
elements of $G$ and $g \in G$ we 
define $ gTg^{-1}:=[gg_1g^{-1},\ldots,gg_rg^{-1}]$.

If
$A=[m_1,\ldots,m_r]\in\NN^r$ is an $r$-tuple of natural numbers with
$2\le m_1\le\ldots \le m_r$ then the spherical system of generators
$T=[g_1,\ldots,g_r]$ is said to have {\it type} 
$A=[m_1,\ldots,m_r]$ if there is a permutation $\tau\in {\mathfrak S}_r$ such
that 
$${\rm ord}(g_1)=m_{\tau(1)},\ldots, {\rm ord}(g_r)=m_{\tau(r)}$$
holds. Here ${\rm ord}(g)$ is the order of the element $g\in G$.
The spherical system of generators $T = [g_1,\ldots,g_r]$ is said to be {\em
  ordered} if $2 \leq {\rm ord}(g_1) \leq \ldots \leq {\rm ord(g_r)}$.

Given a spherical system of generators $T=[g_1,\ldots,g_r]$ of $G$ we define
\begin{equation}
\Sigma(T)=\Sigma([g_1,\ldots,g_r]):=\bigcup_{g\in G}\, \bigcup_{j=0}^\infty\,
\bigcup_{i=1}^r \ \{\, g\cdot g_i^j\cdot g^{-1}\}
\end{equation}
to be the union of all conjugates of the 
cyclic subgroups generated by the elements
$g_1,\ldots, g_r$.
A pair of spherical systems of generators ($T_1,T_2$) of $G$ is called
{\it disjoint} if
$$\Sigma(T_1)\cap \Sigma(T_2)=\{\, 1\,\}.$$ 

\begin{definition}
Consider  a $r$-tuple $A_1=[m_{(1,1)},\ldots,m_{(1,r)}]$   and
a $s$-tuple $A_2=[m_{(2,1)},\ldots, m_{(2,s)}]$  of natural numbers with
$2\le m_{(1,1)}\le\ldots \le m_{(1,r)}$ and $2\le m_{(2,1)}\le\ldots \le m_{(2,s)}$.
An {\em unmixed ramification structure of type $(A_1,A_2)$ for $G$} 
is a 
disjoint pair ($T_1,T_2$) of spherical systems of generators of $G$, 
such that $T_1$ has type $A_1$ and $T_2$ has type $A_2$.
We define ${\cal B}(G;A_1,A_2)$ to be the set of unmixed ramification 
structures 
of type $(A_1,A_2)$ for $G$.
\end{definition}

\begin{definition}\label{defimi}
Let $A=[m_1,\ldots,m_r]$ be a $r$-tuple of natural numbers with
$2\le m_1\le\ldots \le m_r$. 
A {\em mixed ramification structure of type $A$ for $G$} is a pair
$(H,T)$ where $H$ is a subgroup of index $2$ in $G$ and $T=[g_1,\ldots,g_r]$ 
is a $r$-tuple of elements
of $G$ such that the following hold
\begin{itemize}
\item $T$ is a spherical system of generators of $H$ of type $A$,
\item for every $g\in G\setminus H$, the $r$-tuples 
$T$ and $gTg^{-1}=[gg_1g^{-1},\ldots,gg_rg^{-1}]$ are disjoint,
\item for every $g\in G\setminus H$ we have $g^2\notin\Sigma(T)$. 
\end{itemize}
We define ${\cal B}(G;A)$ to be the set of mixed ramification structures 
of type $A$ for $G$.
\end{definition}
We shall now establish certain equivalence relations on the sets 
${\cal B}(G;A_1,A_2)$  and ${\cal B}(G;A)$ of  ramification
structures of a finite group $G$, which  reflect the deformation
equivalence of the surfaces admitting such ramification structures.
 This equivalence relation will be
used in section \ref{modu}.

Let $r$ be a natural number and consider the braid group
\begin{equation}
{\bf B}_r:=\left\langle\, \sigma_1,\ldots,\sigma_{r-1}\, \Mid
\begin{matrix}
\sigma_i\sigma_j=\sigma_j\sigma_i\ {\rm if}\ |i-j|>1,\cr
\sigma_i\sigma_{i+1}\sigma_i=\sigma_{i+1}\sigma_{i}\sigma_{i+1}
\end{matrix}
\right\rangle .
\end{equation}
 We shall define now an action of ${\bf B}_r$ on the set of
$r$-tuples of elements of $G$. This action corresponds to the standard
embedding of ${\bf B}_r$ into the automorphism group of a free group on $r$
generators. 

Let $T=[g_1,\ldots,g_r]$ be a $r$-tuple of elements of $G$ and $1\le i\le r-1$.
Define $\sigma_i(T)$ by
\begin{equation}
\sigma_i(T):=[g_1,\ldots,g_{i-1},g_i\cdot g_{i+1}\cdot
g_i^{-1},g_i,g_{i+2},\ldots ,g_{r}]
\end{equation}

It is well known and also easy to see that 

i) the braid relations are satisfied,

ii) the group ${\bf B}_r$ maps spherical systems of generators to  
spherical systems of generators, preserving the type. 

Also the automorphism
group ${\rm Aut}(G)$ of $G$ acts on the set of spherical systems of generators
of a fixed type by simultaneous application of an automorphism to the
coordinates of a tuple.

Given $(\gamma_1,\gamma_2,\varphi)\in 
{\bf B}_r\times {\bf B}_s\times {\rm Aut}(G)$ and $(T_1,T_2)\in 
{\cal B}(G;A_1,A_2)$, where $T_1$ has length $r$ and $T_2$ has length $s$,
 we set
\begin{equation}\label{act}
(\gamma_1,\gamma_2,\varphi)\cdot (T_1,T_2):=(\varphi(\gamma_1(T_1)),
\varphi(\gamma_2(T_2))).
\end{equation}
A one moment consideration shows that (\ref{act}) leads to an action  
of ${\bf B}_r\times {\bf B}_s\times {\rm Aut}(G)$ on 
${\cal B}(G;A_1,A_2)$. 

Given $(\gamma,\varphi)\in 
{\bf B}_r\times {\rm Aut}(G)$ and $(H,T)\in 
{\cal B}(G;A)$, where $T$ has length $r$, we set
\begin{equation}\label{act2}
(\gamma,\varphi)\cdot (H,T):=(\varphi(H),\varphi(\gamma(T))).
\end{equation}
Formula (\ref{act2}) leads to an action  
of ${\bf B}_r\times {\rm Aut}(G)$ on 
${\cal B}(G;A)$. 

In Section \ref{basi} we will associate to a surface $S$ 
isogenous to a product of unmixed (resp.mixed)
type with $q=0$ an equivalence class of an unmixed (resp. mixed) ramification
structure for
$G=G(S)$. In Section \ref{basi} we shall also conversely see that an 
unmixed (resp. mixed) ramification structure
for a finite group $G$ gives a surface $S$ 
isogenous to a product of unmixed (resp.mixed)
type with $q=0$. 
The equivalence classes (the orbits of the respective actions) determine also exactly
the irreducible components in the corresponding moduli space. This will be
applied in Section \ref{modu}.

\medskip
In the following sections polygonal groups will play an important role.
We  give their definition right away.
Let $A:=[m_1,\ldots,m_r]$ be a $r$-tuple of natural numbers $\geq 2$. The
polygonal group
$\BT(m_1,\ldots,m_r)$ is defined by generators and relations as 
\begin{equation}
\BT(m_1,\ldots,m_r):=\langle\, t_1,\ldots,t_r\Mid t_1t_2\ldots
t_r=1=t_1^{m_1}=\ldots t_r^{m_r}\, \rangle . 
\end{equation}
These groups are important for us since every finite group which has a
spherical system of generators of type $A$ is (in the obvious way) a quotient
group of $\BT(m_1,\ldots,m_r)$.

%%%%%%%%%%%%%%%%%%%%%%%%%%%%%%%%%%%%%%%%%%%%%%%%%%%%%%%%%
\subsection{Tuples}\label{tup}
%%%%%%%%%%%%%%%%%%%%%%%%%%%%%%%%%%%%%%%%%%%%%%%%%%%%%%%%%

In this section we classify $r$-tuples of natural numbers satisfying certain
arithmetic conditions. The lists of these tuples will be of importance in our
later classification program of surfaces.

For an $r$-tuple ($r\in\NN$) $[m_1,\ldots,m_r]\in \NN^r$ define
the orbifold canonical degree as
\begin{equation}
\Theta([m_1,\ldots,m_r]):=-2+\sum_{i=1}^r\left(1-\frac{1}{m_i}\right)
\end{equation}
In the following we define properties of tuples of natural numbers which are
satisfied by the tuples of orders of  the spherical systems of generators 
occurring in the unmixed case. 
\begin{definition}

I) Given $A=[m_1,\ldots,m_r]\in {\cal N}_r$  define
\begin{equation}
\alpha([m_1,\ldots,m_r]):=\frac{2}{\Theta([m_1,\ldots,m_r])}=
\frac{2}{-2+\sum_{i=1}^r\left(1-\frac{1}{m_i}\right)}
\end{equation}
and use the notation $[m_1,\ldots,m_r]_{\alpha(A)}$.

II) For $r\in\BN$ with $r\ge 3$ let ${\cal N}_r$ be the set of $r$-tuples
$[m_1,\ldots,m_r]\in \BN^r$ which satisfy:
\begin{itemize}
\item[{\rm (i):}] $2\le m_1\le\ldots \le m_r$,
\item[{\rm (ii):}] $\Theta([m_1,\ldots,m_r]) > 0$,
\item[{\rm (iii):}]
$\alpha([m_1,\ldots,m_r]):=\frac{2}{\Theta([m_1,\ldots,m_r])}\in \BN$.
\item[{\rm (iv):}] $m_r\le
\alpha([m_1,\ldots,m_r]):=\frac{2}{\Theta([m_1,\ldots,m_r])}$,
\end{itemize}
We set ${\cal N}:=\cup_{i=3}^\infty\, {\cal N}_i$.
\end{definition}

We shall now give a simple classification result for the tuples in ${\cal N}$.
\begin{prop}\label{proputup}
We have ${\cal N}_r=\emptyset$ for $r\ge 7$. The sets ${\cal N}_3$, 
${\cal N}_4$, ${\cal N}_5$, ${\cal N}_6$, are finite and 
\begin{equation} 
{\cal N}_3=\left\{
\begin{matrix}
[ 2 ,  3 , 7 ]_{ 84 }, & [ 2 ,  3 , 8 ]_{ 48 }, & [ 2 ,  4 , 5 ]_{ 40 }, &
[ 2 ,  3 , 9 ]_{ 36 }, & [ 2 ,  3 , 10 ]_{ 30 }, \cr
[ 2 ,  3 , 12 ]_{ 24 }, & [ 2 ,  4 , 6 ]_{ 24 }, & [ 3 ,  3 , 4 ]_{ 24 }, &
[ 2 ,  3 , 14 ]_{ 21 }, & [ 2 ,  3 , 15 ]_{ 20 }, \cr
[ 2 ,  5 , 5 ]_{ 20 }, & [ 2 ,  3 , 18 ]_{ 18 }, & [ 2 ,  4 , 8 ]_{ 16 }, &
[ 2 ,  5 , 6 ]_{ 15 }, & [ 3 ,  3 , 5 ]_{ 15 }, \cr 
[ 2 ,  4 , 12 ]_{ 12 }, & [ 2 ,  6 , 6 ]_{ 12 }, & [ 3 ,  3 , 6 ]_{ 12 }, &
[ 3 ,  4 , 4 ]_{ 12 }, & [ 2 ,  5 , 10 ]_{ 10 }, \cr
[ 2 ,  6 , 9 ]_{ 9 }, & [ 3 ,  3 , 9 ]_{ 9 }, & [ 2 ,  8 , 8 ]_{ 8 }, &
[ 3 ,  4 , 6 ]_{ 8 }, & [ 4 ,  4 , 4 ]_{ 8 },\cr
[ 3 ,  6 , 6 ]_{ 6 }, & [ 4 ,  4 , 6 ]_{ 6 }, & [ 5 ,  5 , 5 ]_{ 5 }
\end{matrix} 
\right\},
\end{equation}
\begin{equation}
{\cal N}_4=\left\{
\begin{matrix}
[ 2,2,2 , 3 ]_{ 12 }, & [ 2 ,  2 , 2 , 4 ]_{ 8 }, & 
[ 2 ,  2 , 2 , 6 ]_{ 6 }, & [ 2 ,  2 , 3 , 3 ]_{ 6 }, \cr
[ 2 ,  2 , 4 , 4 ]_{ 4 }, & [ 2 ,  3 , 3 , 3 ]_{ 4 }, &
[ 3 ,  3 , 3 , 3 ]_{ 3 } &
\end{matrix} 
\right\},
\end{equation}
\begin{equation}
{\cal N}_5=\left\{
[ 2 ,2 , 2 , 2 , 2 ]_{4},\ 
[ 2 ,2 , 2 , 2 , 3 ]_{3}
\right\},\qquad
{\cal N}_6=\left\{
[ 2 ,2 , 2,2 , 2,2]_{2}
\right\}.
\end{equation}
\end{prop}
\Proof
Suppose that $[m_1,\ldots,m_r]$ is a tuple of natural numbers in ${\cal N}_r$.
From condition (iv) we get
\begin{equation}\label{we1}
\sum_{i=1}^r\, (1- \frac{1}{m_i}) \le  2 + \frac{2}{m_r} \leq 3.
\end{equation}
Using $2\le m_i$ for $i=1,\ldots, r$ we obtain $r\le 6$ and 
$ r=6 \Leftrightarrow m_i = 2 \ \forall i$. In particular
 ${\cal N}_r$ is empty for $r\ge 7$.

Let us treat the case $r=3$ next. In this case we have $m_2\ge 3$ since
otherwise $\Theta([m_1,m_2,m_3])$ is negative which contradicts condition
(ii). An application of (\ref{we1}) using $m_1\ge 2$ gives $m_3\le 18$. 
By a quick computer search through the remaining tuples (or just by hand)
we obtain
the finite set ${\cal N}_3$. 
 
In the cases $r=3,\, 4,\, 5,\, 6$ we infer from (\ref{we1}) that 
$\sum_{i=1}^{r-2}\, (1- \frac{1}{m_i}) \le  \frac{1}{m_{r-1}} + \frac{3}{m_r}$,
whence $\frac{r-2}{2}\le \frac{4}{m_{r-1}}$ and $\frac{r-3}{2}\le \frac{3}{m_{r}}$.
These inequalities imply $ m_{r-1} \leq  \frac{8}{r-2} \leq 4$ 
and
$m_r \leq  \frac{6}{r-3} \leq 6$. The remaining tuples can again 
be quickly searched
by computer or by hand to obtain the above lists for
${\cal N}_4$, ${\cal N}_5$, ${\cal N}_6$.
\qed

\medskip

In the following we define properties of a tuple of natural numbers which are
satisfied by the tuple of orders of  the spherical system of generators 
occurring in the mixed case. 

\begin{definition}
I) Given $A=[m_1,\ldots,m_r]\in {\cal M}_r$  define
\begin{equation}
\beta([m_1,\ldots,m_r]):=\frac{4}{\Theta([m_1,\ldots,m_r])}=
\frac{4}{-2+\sum_{i=1}^r\left(1-\frac{1}{m_i}\right)}
\end{equation}
and use the notation $[m_1,\ldots,m_r]_{\beta(A)}$.

II) For $r\in\BN$ with $r\ge 3$ let ${\cal M}_r$ be the set of $r$-tuples
$[m_1,\ldots,m_r]\in \BN^r$ which satisfy:
\begin{itemize}
\item[{\rm (i):}] $2\le m_1\le\ldots \le m_r$,
\item[{\rm (ii):}] $\Theta([m_1,\ldots,m_r]) > 0$,
\item[{\rm (iii):}] $m_r\le
\beta([m_1,\ldots,m_r])=\frac{4}{\Theta([m_1,\ldots,m_r])}$,
\item[{\rm (iv):}] $\beta([m_1,\ldots,m_r])=\frac{4}{\Theta([m_1,\ldots,m_r])}\in
\BN$.
\item[{\rm (v):}] $\beta([m_1,\ldots,m_r])$ is even and 
$\beta([m_1,\ldots,m_r])^2/2$ is divisible by $m_i$ for $i=1,\ldots,m_r$.
\end{itemize}
Define further ${\cal M}:=\cup_{i=3}^\infty\, {\cal M}_i$.
\end{definition}

We shall now give a simple classification result for the tuples in ${\cal M}$.
\begin{prop}\label{propmtup}
We have ${\cal M}_r=\emptyset$ for $r\ge 9$ or $r=7$.
The sets ${\cal M}_3$, 
${\cal M}_4$, ${\cal M}_5$, ${\cal M}_6,$, ${\cal M}_8$, 
are finite and 
\begin{equation}
{\cal M}_3=\left\{
\begin{matrix}
[2,3,7]_{168}, &  [2,3,8]_{96}, & [2,4,5]_{80}, & [2,3,9]_{72}, &  
[2,3,10]_{60}, \cr 
[ 2 , 3 ,  12 ]_{48}, & [ 2 , 4 ,  6]_{ 48 }, & [ 3 , 3 ,  4]_{ 48 }, &
[ 2 , 3 ,  14]_{ 42 }, &  [ 2 , 5 ,  5]_{40},\cr 
[ 2 , 3 ,  18]_{ 36 }, & [ 2 , 4 ,  8]_{ 32 }, & [ 2 , 3 ,  30]_{30 }, & 
[ 2 , 5 ,  6 ]_{ 30 }, & [ 3 , 3 ,  5 ]_{ 30 },\cr  
[ 2 , 4 ,  12]_{ 24 }, & [ 2 , 6 ,  6]_{ 24 } & [ 3 , 3 ,  6]_{ 24 }, & 
[ 3 , 4 ,  4 ]_{ 24 }, & [ 2 , 4 ,  20]_{ 20 },\cr 
[ 2 , 5 ,  10]_{ 20 }, & [ 2 , 6 ,  9]_{ 18 }, & [ 3 , 3 ,  9]_{ 18 },  & 
[ 2 , 8 ,  8 ]_{ 16 }, & [ 4 , 4 ,  4]_{ 16 } \cr
 [ 2 , 7 ,  14]_{ 14 }, & [ 2 , 12 ,  12]_{ 12 }, & [ 3 ,4 ,12]_{ 12 }, &
[ 3 , 6 ,  6]_{ 12 }, & [ 4 , 4 ,  6 ]_{ 12 }, \cr 
[ 5 , 5 ,  5]_{ 10 }, &  [ 4 , 8 ,  8]_{ 8 } &  & & 
\end{matrix} 
\right\},
\end{equation}
\begin{equation}
{\cal M}_4=\left\{
\begin{matrix}
[  2 ,  2 ,  2 ,  3]_{24}, & [  2 ,  2 ,  2 ,  4 ]_{   16 }, &
[  2 ,  2 ,  2 ,  6 ]_{   12 }, & [  2 ,  2 ,  3 ,  3 ]_{  12 }, \cr
[  2 ,  2 ,  2 ,  10 ]_{  10 }, & [  2 ,  2 ,  4 ,  4 ]_{   8 }, &  
[  2 ,  2 ,  6 ,  6 ]_{   6 }, &  [  3 ,  3 ,  3 ,  3 ]_{   6 },\cr 
[  2 ,  3 ,  3 ,  6 ]_{   6 } & [  4 ,  4 ,  4 ,  4 ]_{   4 } & & 
\end{matrix} 
\right\},
\end{equation}
\begin{equation}
{\cal M}_5=\left\{
\begin{matrix}
[ 2 , 2 , 2 , 2 , 2 ]_{ 8 }, &
[ 2 , 2 , 2 , 2 , 3 ]_{ 6 }, &
[ 2 , 2 , 2 , 4 , 4 ]_{ 4 }
\end{matrix} 
\right\},
\end{equation}
\begin{equation}
{\cal M}_6=\left\{
[2,    2,    2,    2,    2,    2]_{4}
\right\},
\quad {\cal M}_8=\left\{[2,2,2,2,2,2,2,2]_2\right\}
\end{equation}
\end{prop}
We skip the proof since it is similar to that of Proposition \ref{proputup}.

%%%%%%%%%%%%%%%%%%%%%%%%%%%%%%%%%%%%%%%%%%%%%%%%%%%%%%%
\section{Basics on surfaces isogenous to a product}\label{basi}
%%%%%%%%%%%%%%%%%%%%%%%%%%%%%%%%%%%%%%%%%%%%%%%%%%%%%%%
Throughout this section we assume that $S$ is a surface of general type.
We recall first the notion of surfaces isogenous to a product of 
curves. By Proposition 3.11 of
\cite{cat00} the following two properties for a surface of general type are 
equivalent.

\begin{definition}\label{defbasi}  
A surface $S$ of general type is said to be {\em isogenous to a 
product} if and only if one of the following two equivalent 
conditions is satisfied.
\begin{itemize}
\item $S$ admits a finite unramified covering which is isomorphic to a 
product of curves (of genera at  least two),  
\item $S$ is a quotient  $S := (C_1 \times C_2) / G$, 
where the $C_i$'s are
curves of genus at least two, and $G$ is a finite group acting freely on
$C_1 \times C_2$.
\end{itemize}
\end{definition}

It is shown in \cite{cat00} that every such 
surface isogenous to a product has 
a unique minimal realization $S := (C_1 \times C_2) / G$ 
(i.e., the genera $g(C_1),\, g(C_2)$
of the two curves  $C_1 , C_2$ are minimal). 

It can further be shown (see \cite{cat00}) that the action 
of $G$ on $C_1 \times C_2$ in the second condition of the above
definition respects the product decomposition, i.e., the elements of $G$  
either interchange the factors or act independently on both factors. 
\begin{definition}
Let $S$ be a surface isogenous to a 
product with minimal realisation $S=(C_1\times C_2)/G$.
We say that $S$ is a {\rm  mixed case} if the action of $G$ 
exchanges the two factors (and
then $C_1 , C_2$ are isomorphic), and an {\rm   unmixed case} if
$G$ acts via a diagonal action.
\end{definition}
We shall associate to a surface $S$ now certain algebraic data. This approach
is taken from \cite{cat00} where a much more detailed discussion can be found. 
We first take a minimal realisation of $S$ as $S=(C_1\times C_2)/G$
and define
\begin{equation}
G(S):=G.
\end{equation}
Suppose we are in the unmixed case. Then $q(S)=0$ implies that 
$C_1/G(S)=C_2/G(S)=\PP^1$, i.e., we have two ramified Galois coverings 
\begin{equation}
C_1\to \PP^1, \qquad C_2\to \PP^1
\end{equation}
with Galois group $G$ (see \cite{miranda}, Section 4 for explanations).  
Let $\{ P_1,\ldots, P_r\}\subset \PP^1$ be the set 
of branch points of the first covering. Choose a base point $P$
in $\PP^1$ distinct from them. 
Choose a 
geometric basis  $\gamma_1, \ldots
\gamma_r$ of $\pi_1( \PP^1 - \{P_1, \dots P_r \}) $ ($\gamma_i$ is a 
simple counterclockwise loop
around
$P_i$, and they follow each other by counterclockwise ordering around 
the base point). 
Notice that $\gamma_1\cdot\ldots \cdot\gamma_r=1$.
Choose a monodromy representation, i.e., a surjective homomorphism 
$$\psi: \pi_1( \PP^1 - \{P_1, \dots P_n \})  \to G.$$ 
Notice that only the kernel of $\psi$ is uniquely determined
by the covering.
Then the elements $\psi(\gamma_1),\ldots,\psi(\gamma_r)$ form a 
spherical system of generators of $G$.  

Now, the mapping class group of the sphere $\pi_0( Diff ( \PP^1 - \{P_1, 
\dots P_n \}))$, which is a quotient of the braid group ${\bf B}_n$, operates
on such homomorphisms, and their orbits are called Hurwitz 
equivalence classes of spherical
systems of generators. This action is the one which was already 
described in the previous section. We use this action in
order to assume without loss of generality that $T_1(S):=
[\psi(\gamma_1),\ldots,\psi(\gamma_r)]$ is an ordered spherical system of
generators.

We apply the same principle to the second covering and obtain another ordered
spherical system of generators $T_2(S)$ of $G$.

Since the action of $G$ on $C_1 \times C_2$ is free we have
$\Sigma(T_1(S)) \cap\Sigma(T_2(S)) = \{\, 1\,\},$ i.e., the two systems are
disjoint.

\medskip
Let $S$ be a surface isogenous to a product, of unmixed type and with $q(S)=0$.
Then we have attached to $S$ its finite group $G=G(S)$ (up to isomorphism) and 
a pair ${\cal T}(S)=(T_1(S),T_2(S))\in {\cal B}(G;A_1(S),A_2(S))$ of 
an uniquely defined ordered 
type $(A_1(S),A_2(S))$. 

We show now that the tuples $T_1(S)$, $T_2(S)$ attached to a surface $S$
isogenous to a product, of unmixed type and with $p_g(S)=0$ satisfy the
properties of section \ref{tup}, i.e., that 
they are contained in $\mathcal{N}$.
\begin{prop}\label{unmitup}
Let $S$ be a surface isogenous to a product, of 
unmixed type and with $p_g(S)=0$. Let $A_1(S)=[m_1,\ldots,m_r]$, 
$A_2(S)=[n_1,\ldots,n_s]$ be the two ordered types attached to $S$ as above.
We have 
\begin{itemize}
\item $\Theta(A_1(S)),\, \Theta(A_2(S)) > 0$,
\item $m_r\le \frac{2}{\Theta(A_1(S))}$, $n_s\le \frac{2}{\Theta(A_2(S))}$,
\item $\frac{2}{\Theta(A_1(S))},\, \frac{2}{\Theta(A_2(S))}   \in \BN$.
\end{itemize}
\end{prop}
{\it Proof.} Since $S$ is isogenous to a product we can represent $S$ as
$$S=(C_1\times C_2)/G(S)$$
where $C_1,\, C_2$ are two smooth projective curves with genera 
$g(C_1),\, g(C_2)\ge
2$ where the finite group $G(S)$ acts without  fixed points
and via an action preserving the product on $C_1\times C_2$.
Since $q(S)=0$ we have $C_1/G(S)\cong  C_2/G(S)\cong \mathbb{P}^1$ and 
the Hurwitz formula implies
\begin{equation}\label{e1} 
|G|\left(-2 + \sum_j^r \left(1 - \frac{1}{m_j}\right)\right)=2(g(C_1) - 1),
\end{equation}
\begin{equation}\label{e2}
|G|\left(-2 + \sum_j^s \left(1 - \frac{1}{n_j}\right)\right) =2(g(C_2) - 1)
\end{equation} 
This establishes $\Theta(A_1(S)),\, \Theta(A_2(S)) >  0$, because 
$g(C_1),\, g(C_2)\ge 2$.

We have
$$ K_S^2 = \frac{K^2_{C_1\times C_2}}{|G|} = \frac{8
  (g(C_1) - 1)(g(C_2) -1)}{|G|} = 8 \chi (\mathcal{O}_S)= 8,
$$
where the last equality holds since $p_g = 0$.

Therefore
$$|G(S)|=(g(C_1)-1)(g(C_2)-1)$$ and 

using formulas (\ref{e1}), (\ref{e2})  we get
$$
\frac{2}{\Theta(A_1)}=g(C_2)-1,\ \frac{2}{\Theta(A_2)}=g(C_1)-1  
\in \mathbb{N}.
$$
This establishes the third property of $A_1$, $A_2$.

To prove the second property assume that 
$$m_r> \frac{2}{\Theta(A_1(S))}=g(C_2)-1.$$
If $T(S)=[g_1,\ldots,g_r]$ then $g_r$ has order $m_r$ and 
we know that  it acts with a fixed
point on $C_1$. Hence the cyclic group $\langle g_r\rangle$ 
should have no fixed points
on $C_2$. 
Let $C := C_2 / \langle g_r\rangle $ be the quotient.
By Hurwitz' formula we get:
$$
2g(C) - 2 = \frac{2g(C_2) - 2}{m_r} < 2.
$$
Therefore $g(C) \in \{ 0, 1 \}$, which  contradicts 
the freeness of the  action of $\langle g_r\rangle$ on $C_2$ (recall that $\BP^1$
has no unramified coverings and an unramified covering of an elliptic curve is again
an elliptic curve, while $C_2$ has genus $\geq 2$). 

\qed 

\medskip
Let $S$ be a surface isogenous to a product, of mixed type and with $q(S)=0$.
Then we can attach to $S$ its finite group $G=G(S)$ and 
a pair ${\cal T}(S)=(H(S),T(S))\in {\cal B}(G;A)$ of 
an uniquely defined ordered
type $A(S)$.

Here we get the following

\begin{prop}
Let $S$ be a surface isogenous to a product, of mixed type and with
$p_g(S)=0$. Let $A(S)=[m_1,\ldots,m_r]$ be the ordered type attached to $S$.
We have 
\begin{itemize}
\item $\Theta(A(S)) \ne 0$,
\item $m_r\le \frac{4}{\Theta(A(S))}$,
\item $\beta(A(S)):= \frac{4}{\Theta(A(S))}  \in \BN$,
\item $\beta(A(S))$ is even and 
$\beta(A(S))^2/2$ is divisible by $m_i$ for $i=1,\ldots,m_r$.
\end{itemize}
\end{prop}

\Proof
Noting that $H(S)$ has an unmixed ramification structure of type
$(A,A)$ yielding a surface with invariants $K_S^2= 16 = 8 \chi$, the first
three properties are proven in the same way as in proposition
\ref{unmitup}. For the last property observe that $G$ has order
$\beta(A(S))^2$ and has a subgroup of index $2$. Moreover, $|H(S)| =
\beta(S(S))^2/2$ and has a spherical system of generators of type $A(S) = 
[m_1,\ldots,m_r]$.
\qed

\medskip
So far we have discussed the ramification structure associated to a surface
isogenous to a product. There is also a way back from ramification structures
to surfaces. This construction relies on the Riemann existence theorem (see
\cite{cat00} for  details). More precisely we have 
\begin{prop}
Let $G$ be a finite group.
let $A_1=[m_{11},\ldots,m_{1r}]$ be a $r$-tuple and
$A_2=[m_{21},\ldots,m_{2s}]$ a $s$-tuple of natural numbers with
$2\le m_{11}\le\ldots \le m_{1r}$ and $2\le m_{21}\le\ldots \le m_{2s}$.
Then for any ramification structure ${\cal T}\in {\cal B}(G;A_1,A_2)$ 
there is a
surface isogenous to a product with $G(S)=G$ and ${\cal T}(S)={\cal T}$.
\end{prop}

An analogous existence result holds in the mixed case also.

%%%%%%%%%%%%%%%%%%%%%%%%%%%%%%%%%%%%%%%%%%%%%
\section{The unmixed case, classification of the groups}\label{classiumi}
%%%%%%%%%%%%%%%%%%%%%%%%%%%%%%%%%%%%%%%%%%%%%%

This section is devoted to the classification of all finite groups $G$ 
admitting
an unmixed ramification structure of type $(A,B)$ with
$A,\, B\in {\cal N}$. The result is summarized in the following

\begin{prop}
The only finite groups $G$ admitting
an unmixed ramification structure of type $(A,B)$ with
$A,\, B\in {\cal N}$ are those in the following table:
\begin{center}
\begin{tabular}{|c|c|c|c|}
\hline
$G$  & $|G|$ & $A$ & $B$\\
\hline
${\mathfrak A}_5$ & $60$ & $[ 2, 5, 5]_{20}$ & $[ 3, 3, 3, 3]_{3}$\\
${\mathfrak A}_5$ & $60$ & $[ 5, 5, 5]_{ 5 }$ & $[ 2, 2, 2, 3]_{12}$ \\
${\mathfrak A}_5$ & $60$ & $[ 3, 3, 5]_{15}$ & $[ 2, 2, 2, 2]_{4}$ \\
${\mathfrak S}_4\times \Z_2$ & $48$ & $[ 2, 4, 6]_{24}$ & $[ 2, 2, 2, 2, 2, 2]_2$ \\
G($32$) & $32$ & $[ 2, 2, 4, 4]_4$ & $[ 2, 2, 2, 4]_8$ \\
$\Z_5^2$ & $25$ & $[5,5,5]_5$ & $[5,5,5]_5$ \\
${\mathfrak S}_4$ & $24$ & $[ 3, 4, 4]_{12}$ & $[ 2, 2, 2, 2, 2, 2]_2$ \\
G($16$) & $16$ & $[ 2, 2, 4, 4]_4$ & $[ 2, 2, 4, 4]_4$ \\
${\rm D}_4\times\Z_2$ & $16$ & $[ 2, 2, 2, 4]_8$ & $[ 2, 2, 2, 2, 2, 2]_2$\\
$\Z_2^4$ & $16$ & $[ 2, 2, 2, 2, 2]_4$ & $[ 2, 2, 2, 2, 2]_4$ \\
$\Z_3^2$ & $9$ & $[ 3, 3, 3, 3]_3$ & $[ 3,3, 3, 3]_3$ \\ 
$\Z_2^3$ & $8$ & $[ 2, 2, 2, 2, 2]_4$ & $[ 2, 2, 2, 2, 2, 2]_2$ \\ 
\hline
\end{tabular}
\end{center}
\end{prop}

\medskip
The proof relies heavily on the use of the MAGMA-library containing either
permutational representations or polycyclic presentations of all groups of order less than
$2000$ with the exception of the order $1024$. We proceed as follows. 
We consider each type $(A,B)$ separately going through
the finite list of Proposition \ref{proputup}. Assume that there is a group 
$G$ admitting
an unmixed ramification structure of type $(A,B)$:
then $|G|=\alpha(A)\alpha(B)$.
If this order is less than $2000$ we just go through the MAGMA-library
and search for groups which have a disjoint pair of systems of spherical generators
of type $(A,B)$. There is a huge number of groups to check but there are methods
to speed up the computation. 
These will be described in the next two subsections
where we also exhibit the arguments for $|G|=\alpha(A)\alpha(B)>2000$.
Sometimes we shall have to talk about individual groups in the
MAGMA-library. Here we use the terminology of MAGMA, i.e.,
$\SmG(a,b)$ denotes the group of order $a$ having number $b$ in the list.

A simple but useful observation is 
\begin{lemma}
Let $(T_1, T_2)$ be a 
disjoint pair of spherical systems of generators of a finite group $G$.
Then, for every $g\in G$, $gT_1g^{-1}$, $T_2$ and 
$T_1$, $gT_2g^{-1}$ are also  
disjoint pairs of spherical systems of generators of $G$.
\end{lemma}
\Proof $\Sigma(T_1)$ and $\Sigma(T_2)$ are unions of conjugacy classes. 
\qed

\begin{rem}
We will often use without explicit mention the
$p^{\alpha}q^{\beta}$-Theorem (of Burnside) saying that every group of order  
$p^{\alpha}q^{\beta}$ ($p,\, q$ primes) is solvable (cf. \cite{burnart} and
\cite{burnlibro}, p.323). 
\end{rem}

%%%%%%%%%%%%%%%%%%%%%%%%%%%%%%%%%%%%%%%%
\subsection{The case: $A=[2,3,7]_{84},\ B\in{\cal N}$}
%%%%%%%%%%%%%%%%%%%%%%%%%%%%%%%%%%%%%%%%

In this section we consider the case $A=[2,3,7]_{84}$ and prove
\begin{prop}\label{prop84}
There is no finite group $G$ having an unmixed ramification structure
of type $(A,B)$ with  $A=[2,3,7]_{84}$ and $B\in{\cal N}$ arbitrary.
\end{prop}

Since the finite group $G$ satisfying the conditions of Proposition
\ref{prop84} has a system of generators of type $[2,3,7]$ it has to be a
non-trivial perfect group. 
Recall that a group $G$ is called 
perfect if $G=G'$ where $G'$ is the commutator
subgroup of $G$. Notice that every quotient group of a perfect group is again
perfect. Running through the relevant MAGMA library it can quickly be checked
whether there is a perfect group of some cardinality. We shall often exploit
\begin{rem}\label{simplequot}
Let $G$ be a non trivial finite perfect group. 
Then $G$ has a non abelian simple group
$Q$ as quotient.
\end{rem}
We shall also make use of 
\begin{comp}\label{ord}
The only non abelian simple groups whose order divides $2016$, 
$3024$, $4032$, $7056$ are $\PSL(2,\BF_7)$ and $\PSL(2,\BF_8)$.
\end{comp}
This can be seen by applying the command 
``SimpleGroupsWithOrderDividing'' of MAGMA to the above numbers. 

Now we shall run through all $B\in{\cal N}$ and 
explain the computations that are performed. 

%%%%%%%%%%%%%%%%%%%%%%%%%%%%%%%%%%%%%%%%%%%%%%%%%%%%%
\subsubsection{$A=[2,3,7]_{84},\ B\in{\cal N},\ \alpha(B)\le 21$}
\label{sususe1}
%%%%%%%%%%%%%%%%%%%%%%%%%%%%%%%%%%%%%%%%%%%%%%%%%%%%%%

The necessary computations can be speeded up enormously by first establishing
the following
\begin{comp}
The only perfect groups of order $84k$, where 
$k\le 21$ is one of the numbers $\alpha(B)$ for $B\in{\cal N}$,
are

\noindent
1. $\SmG(168,42)$ for $k=2$,

\noindent
2. $\SmG(336,114)$ for $k=4$,

\noindent
3. $\SmG(504,156)$ for $k=6$,

\noindent 
4. $\SmG(1344,814)$ for $k=16$,

\noindent 
5. $\SmG(1344,1186)$ for $k=16$.

\noindent
The first four of them have only one conjugacy class of elements of order $2$.
\end{comp}
This can be seen quickly by running through the relevant MAGMA-libraries
checking the IsPerfect-predicate.

We then exclude the first four cases since we verify that for each $B \in
\mathcal{N}$ with $\alpha(B) = 2,4,6,16$ we 
have that $2$ divides one of the $m_i$'s.
We are left with the case $B=[2,4,8]_{16}$ 
(since $\alpha(B) = 16$ iff $B=[2,4,8]_{16}$) and $G=\SmG(1344,1186)$.
By Computational Fact \ref{ord} this 
group is an extension of $\PSL(2,\BF_7)$ by
a kernel of order $8$ and it can be excluded noting that 
$\PSL(2,\BF_7)$ has no spherical system of generators of types 
$[2,4,8]$, $[2,4,4]$, $[2,2,4]$ or $[2,2,2]$.

%%%%%%%%%%%%%%%%%%%%%%%%%%%%%%%%%%%%%%%%%%%%%%%%%
\subsubsection{$A=[2,3,7]_{84},\ B\in{\cal N}_3,\ \alpha(B)=24$}
%%%%%%%%%%%%%%%%%%%%%%%%%%%%%%%%%%%%%%%%%%%%%%%%%

Here the order of $G$ is $2016$. By \ref{simplequot} there is a simple non
abelian quotient $Q$ of $G$, 
which by Computational fact \ref{ord} is isomorphic to $\PSL(2,\BF_7)$ or to
$\PSL(2,\BF_8)$.

Suppose that $Q$ is $\PSL(2,\BF_7)$: 
then the kernel $K$ of the quotient
homomorphism has order $12$. Each group of order $12$ has either a normal
Sylow-$2$-subgroup or a normal Sylow-$3$-subgroup. This Sylow-subgroup, which
we denote by $S$, is characteristic in $K$, hence normal in $G$. 
Suppose that $|S|=3$. Then $G/S$ is perfect and has order
$2016/3=672$, but there are no perfect
groups of this order (checked by MAGMA). If $|S|=4$, then 
$G/S$ has order $2016/4=504$. The only perfect group of order $504$ is 
$\PSL(2,\BF_8)$ which is simple. A contradiction (since $K/S$ is a normal
subgroup of $G/S$).

Suppose that $Q$ is $\PSL(2,\BF_8)$: then the kernel $K$ of the quotient
homomorphism has order $4$. Since $Q$ has to act trivially (by conjugation)
on $K$. Let $K_1$ be a subgroup of order $2$ in $K$. Clearly $K_1$ is normal
in $G$ and $G/K_1$ is perfect of order $1008$. There are no such groups.

This shows that there is no group $G$ of order $2016$ with a spherical 
system of
generators of type $[2,3,7]$.

%%%%%%%%%%%%%%%%%%%%%%%%%%%%%%%%%%%%%%%%%%%%%%%%%
\subsubsection{$A=[2,3,7]_{84},\ B\in{\cal N}_3,\ \alpha(B)=30$}
%%%%%%%%%%%%%%%%%%%%%%%%%%%%%%%%%%%%%%%%%%%%%%%%%

Here the order of $G$ is $2520$ and $B=[2,3,10]_{30}$. We find
\begin{comp}
The only non abelian simple groups of order dividing $2520$ are
${\mathfrak A}_5$, ${\mathfrak A}_6$, ${\mathfrak A}_7$,
$\PSL(2,\BF_7)$ and $\PSL(2,\BF_8)$.
\end{comp}
Let $Q$ be a simple quotient of $G$. 
The cases $Q={\mathfrak A}_5$ and $Q={\mathfrak A}_6$ cannot occur since these
groups have no element of order $7$. 

$Q={\mathfrak A}_7=G$ cannot occur since  
${\mathfrak A}_7$ has only one conjugacy class of elements of order $2$ and hence 
cannot have a disjoint pair of spherical generators of type $(A,B)$.

Suppose now that $Q=\PSL(2,\BF_7)$ or $Q=\PSL(2,\BF_8)$. 
These groups do not have an element of order
$5$, hence $Q$ would have to be a quotient of $\BT(2,3,2)$ which is a dihedral
group. A contradiction.

%%%%%%%%%%%%%%%%%%%%%%%%%%%%%%%%%%%%%%%%%%%%%%%
\subsubsection{$A=[2,3,7]_{84},\ B\in{\cal N}_3,\ \alpha(B)=36$}
%%%%%%%%%%%%%%%%%%%%%%%%%%%%%%%%%%%%%%%%%%%%%%%%%

Here the order of $G$ is $3024$. By \ref{ord} a 
non abelian simple quotient $Q$ of $G$ can only be $\PSL(2,\BF_7)$ or
$\PSL(2,\BF_8)$.

Suppose that $Q=\PSL(2,\BF_7)$, then the kernel $K$ of 
the quotient homomorphism 
from $G$ to $Q$ has order $18$. The Sylow-$3$-subgroup $S$ of $K$ is normal, 
hence characteristic in $K$. It follows that $S$ is normal in $G$. 

We have 
\begin{comp}
There is only one perfect group of order $336$, namely ${\rm SL}(2,\BF_7)$.
\end{comp}
Hence $G/S$ is isomorphic to ${\rm SL}(2,\BF_7)$. This is a contradiction
because $G/S$ would be a quotient of $\BT(2,3,7)$ and 
${\rm SL}(2,\BF_7)$ has only one element of order $2$ which lies in its center.

Suppose that $Q=\PSL(2,\BF_8)$. The kernel $K$ of the quotient homomorphism 
from $G$ to $Q$ has order $6$. Its Sylow-$3$-subgroup $S$ is normal in $G$.
The quotient group $G/S$ has order $1008$ and is perfect. There is no
such group. 

This shows that there is no group $G$ of order $3024$ with a spherical system
of
generators of type $[2,3,7]$.

%%%%%%%%%%%%%%%%%%%%%%%%%%%%%%%%%%%%%%%%%%%%%%%
\subsubsection{$A=[2,3,7]_{84},\ B\in{\cal N}_3,\ \alpha(B)=40$}
%%%%%%%%%%%%%%%%%%%%%%%%%%%%%%%%%%%%%%%%%%%%%%%

Here the order of $G$ is $3360$ and $B=[2,4,5]_{40}$. We find
\begin{comp}
The only non abelian simple groups of order dividing $3360$ are
${\mathfrak A}_5$, and $\PSL(2,\BF_7)$.
\end{comp}
Any minimal non abelian simple quotient $Q$ of $G$ has to be one of these
two groups. Using the fact that 
${\mathfrak A}_5$ has no element of order $7$ and $\PSL(2,\BF_7)$ has no 
element of order $5$ we see that there is no group $G$ in this section.

%%%%%%%%%%%%%%%%%%%%%%%%%%%%%%%%%%%%%%%%%%%%%%%%%%
\subsubsection{$A=[2,3,7]_{84},\ B\in{\cal N}_3,\ \alpha(B)=48$}\label{ss1}
%%%%%%%%%%%%%%%%%%%%%%%%%%%%%%%%%%%%%%%%%%%%%%%%%%

Here the order of $G$ is $4032$ and $B=[2,3,8]_{48}$.
By Computational fact \ref{ord} 
any minimal non abelian simple quotient $Q$ of $G$ has to be one of
$\PSL(2,\BF_7)$ or $\PSL(2,\BF_8)$.
But these two groups do not have a 
spherical system of generators of type $[2,3,8]$, $[2,3,4]$ or $[2,3,2]$. 

%%%%%%%%%%%%%%%%%%%%%%%%%%%%%%%%%%%%%%%%%%%%%%%%%%%%
\subsubsection{$A=[2,3,7]_{84},\ B\in{\cal N}_3,\ \alpha(B)=84$}
%%%%%%%%%%%%%%%%%%%%%%%%%%%%%%%%%%%%%%%%%%%%%%%%%%%%

Here the order of $G$ is $7056$. 
By Computational fact \ref{ord} 
any minimal non abelian simple quotient $Q$ of $G$ has to be
$\PSL(2,\BF_7)$ or $\PSL(2,\BF_8)$.
Let $K$ be the kernel of the
quotient homomorphism. This group has order $42$ or $14$ and 
it has a normal, hence characteristic
Sylow-$7$-subgroup $S$, which has to be normal in $G$. Then $G/S$ has
order $1008$  and is perfect. There are no such groups.

%%%%%%%%%%%%%%%%%%%%%%%%%%%%%%%%%%%%%%%%
\subsection{The case: $A=[2,3,8]_{48},\ B\in{\cal N}$}
%%%%%%%%%%%%%%%%%%%%%%%%%%%%%%%%%%%%%%%%

In this section we treat the case $A=[2,3,8]_{48}$ and prove
\begin{prop}\label{prop48}
There is no finite group $G$ having an unmixed ramification structure
of type $(A,B)$ with  $A=[2,3,8]_{48}$ and $B\in{\cal N}$ arbitrary.
\end{prop}
Let us first consider the case $B\in{\cal N}_3$. All cases except
$B=[2,3,7]_{84}$, $[2,3,8]_{48}$ or $[2,4,8]_{16}$ can be analysed quickly by
MAGMA. In fact, there is no group $G$ 
having simultaneously a spherical system of generators
of type $A$ and one of type $B$ in these cases. The computer calculation is
speeded up enormously by noting that $G$ has to be either perfect or has
abelianisation equal to $\Z_2$.

In the cases $B=[2,3,7]_{84}$, $[2,3,8]_{48}$ the order of $G$ is bigger than
$2000$ and there are no MAGMA-libraries of groups of such order. The
first case was already treated in Section \ref{ss1}, the second is treated below. In case
$B=[2,4,8]_{16}$ the number of groups ($1090235$) to be considered makes the
computation time consuming, so we treat it by a direct argument below.

%%%%%%%%%%%%%%%%%%%%%%%%%%%%%%%%%
\subsubsection{$A=[2,3,8]_{48}$, $B\in {\cal N}_3,\ \alpha(B)=16$}
%%%%%%%%%%%%%%%%%%%%%%%%%%%%%%%%

Here we have $B=[2,4,8]_{16}$. The order of the group $G$ is $768 = 3 \cdot 256$ and
$G$ is solvable.

The group $G$ is a quotient of $\BT(2,3,8)$ and of $\BT(2,4,8)$.
We have 
\begin{equation}\label{tria1}
\BT(2,3,8)^{\rm ab}=\Z_2,\qquad \BT(2,3,8)' \cong \BT(3,3,4).
\end{equation}
The abelianisation of $\BT(3,3,4)$ is $\Z_3$. These facts imply that $G^{\rm
  ab}$ is $\Z_2$ and $(G')^{\rm ab}$ is $\Z_3$, since they are both not perfect.  

The triangle group $\BT(2,4,8)$ has exactly $3$ subgroups of index $2$. They are
isomorphic to $\BT(2,8,8)$, $\BT(4,4,4)$ and $\BT(2,2,2,4)$ and each of them
has a $2$-group as abelianisation. This implies that the 
abelianisation of $G'$ has to be a $2$-group which is a contradiction.

The subgroups of finite index in finitely presented groups (like $\BT(2,4,8)$)
can quickly be analysed by the ``generators and relations'' programs of MAGMA
to obtain results like those just used (alternatively one can use geometric
branching arguments, cf. lemma $(2.3)$ of \cite{cat03}).

%%%%%%%%%%%%%%%%%%%%%%%%%%%%%%%%%%%%%%%%%%%%
\subsubsection{$A=[2,3,8]_{48}$, $B\in {\cal N}_3,\ \alpha(B)=48$}
%%%%%%%%%%%%%%%%%%%%%%%%%%%%%%%%%%%%%%%%%%%%%

In this case we have $B=[2,3,8]_{48}$ and $|G|=2304$.
Every group of this order is solvable. 

We again use (\ref{tria1}) to find that $G$
has to have a subgroup of index 2 with 
abelianisation $\Z_3$. The
Smallgroups-library of MAGMA contains all groups of order 1152. It can
be quickly checked that there is no group of order 1152 with
abelianisation equal to $\Z_3$.

%%%%%%%%%%%%%%%%%%%%%%%%%%%%%%%%%%%%%%%%%%%%
\subsubsection{$A=[2,3,8]_{48}$, $B\in {\cal N}_4,\ {\cal N}_5,\ {\cal N}_6$}
%%%%%%%%%%%%%%%%%%%%%%%%%%%%%%%%%%%%%%%%%%%%%

All these pairs of types $(A,B)$ can be quickly searched by computer to
find
\begin{comp}
There is no finite group $G$ having a pair of sytems of generators
of types $(A,B)$ with  $A=[2,3,8]_{48},\,  B\in {\cal N}_4,$ 
${\cal N}_5$, ${\cal N}_6$ with exception of 
$\SmG(96,64)$ which has non disjoint pairs of spherical systems of 
generators of 
type $(A,B)$ with $B=[2,2,2,2,2,2]_2$.
\end{comp}

%%%%%%%%%%%%%%%%%%%%%%%%%%%%%%%%%%%%%%%%
\subsection{The case: $A,\, B\in{\cal N}_3,\ \alpha(A),\, \alpha(B)\le 40$}
%%%%%%%%%%%%%%%%%%%%%%%%%%%%%%%%%%%%%%%%

All these pairs of types $(A,B)$ can be quickly searched by computer to
find
\begin{comp}\label{prop40}
There is no finite group $G$ having a disjoint pair of spherical 
sytems of generators
of types $(A,B)$ with  $A,\, B\in {\cal N}_3$ and 
$\alpha(A),\, \alpha(B)\le 40$ except for $G=\Z_5^2$ which admits such a
system with $A,\, B=[5,5,5]_5$.
\end{comp}

%%%%%%%%%%%%%%%%%%%%%%%%%%%%%%%%%%%%%%%%
\subsection{The case: $A\in {\cal N}_3,\, \alpha(A)\ne 48,\, 84,\ 
B\in{\cal N}_4,\ {\cal N}_5,\ {\cal N}_6$}
%%%%%%%%%%%%%%%%%%%%%%%%%%%%%%%%%%%%%%%%

All these pairs of types $(A,B)$ can be quickly searched by computer to 
find as only groups $G$ with an unmixed ramification structure the following cases
\begin{itemize}
\item $A=[ 2, 5, 5]_{20},\, B=[ 3, 3, 3, 3]_{ 3}$, $G=\SmG(60,5)$,
\item $A=[ 5, 5, 5]_{5},\, B=[ 2, 2, 2, 3]_{ 3}$, $G=\SmG(60,5)$,
\item $A=[ 3, 3, 5]_{15},\, B=[ 2, 2, 2, 2,2]_{ 4}$, $G=\SmG(60,5)$,
\item $A=[ 2, 4, 6]_{24},\, B=[ 2, 2, 2, 2,2,2]_{ 2}$, $G=\SmG(48,48)$,
\item $A=[ 3, 4, 4]_{12},\, B=[ 2, 2, 2, 2,2,2]_{ 2}$, $G=\SmG(24,12)$.
\end{itemize}
From the description of the groups given by MAGMA it is easy to see that
$\SmG(60,5)$ is the alternating group ${\mathfrak A}_5$, $\SmG(48,48)$ is
${\rm S}_4\times \Z_2$ and $\SmG(24,12)$ is ${\rm S}_4$.

%%%%%%%%%%%%%%%%%%%%%%%%%%%%%%%%%%%%%%%%
\subsection{The case: $A,\, B\in {\cal N}_4,\ {\cal N}_5,\ {\cal N}_6$}
\label{suse44}
%%%%%%%%%%%%%%%%%%%%%%%%%%%%%%%%%%%%%%%%

All these pairs of types $(A,B)$ can be quickly searched by computer to 
find as only groups $G$ with an unmixed ramification structure the following cases
\begin{itemize}
\item $A=[ 2, 2, 2,4]_{8},\, B=[ 2, 2, 4, 4]_{ 4}$, $G=\SmG(32,27)$,
\item $A=[ 2, 2, 4,4]_{4},\, B=[ 2, 2, 4, 4]_{ 4}$, $G=\SmG(16,3)$,
\item $A=[ 2, 2, 2,4]_{8},\, B=[ 2, 2, 2, 2,2,2]_{ 2}$, $G=\SmG(16,11)$,
\item $A=[ 2, 2, 2,2,2]_{4},\, B=[ 2, 2, 2, 2,2]_{ 4}$, $G=\SmG(16,14)$,
\item $A=[ 3, 3, 3,3]_{9},\, B=[ 3, 3, 3, 3]_{ 3}$, $G=\SmG(9,2)$,
\item $A=[ 2, 2, 2,2,2]_{4},\, B=[2,2,2,2,2,2]_{2}$, $G=\SmG(8,5)$.
\end{itemize}
The group $\SmG(8,5)$ is $\Z_2^3$, $\SmG(9,2)$ is $\Z_3^2$, 
$\SmG(16,14)$ is $\Z_2^4$ and $\SmG(16,11)$ is ${\rm D}_4\times \Z_2$ where 
${\rm D}_4$ stands for the dihedral group of order $8$.
A finite presentation of $\SmG(16,3)$ is 
\begin{equation}\label{gr16}
G(16):=\SmG(16,3)=\langle\, g_1,\, g_2,\, g_3,\, g_4\Mid g_1^2=g_4,\, 
g_2^{g_1}=g_2g_3\,\rangle.
\end{equation}

\begin{rem}
The convention here is that the squares of all generators $g_1,\ldots,g_4$
which are not mentioned in the presentation are equal to $1$. If $h_1,\, h_2$
are elements of the group $G$ then $h_1^{h_2}:=h_2^{-1}h_1h_2$. All conjugates
$g_i^{g_j}$ amongst the generators which are not mentioned are equal to $g_i$,
i.e., $g_i$ and $g_j$ commute in this case.
\end{rem} 

A finite presentation of $G(32):=\SmG(32,27)$ is
\begin{equation}\label{gr32}
G(32)=\langle\, g_1,\, g_2,\, g_3,\, g_4,\, g_5\Mid  
g_2^{g_1}=g_2g_4, \, g_3^{g_1}=g_3g_5     \,\rangle.
\end{equation}

%%%%%%%%%%%%%%%%%%%%%%%%%%%%%%%%%%%%%%%%%%%%%%%%%%%%%%%%%
\section{The mixed case, classification of the groups}\label{classimi}
%%%%%%%%%%%%%%%%%%%%%%%%%%%%%%%%%%%%%%%%%%%%%%%%%%%%%%%%%

This section contains the classification of all finite groups $G$ which admit
a mixed ramification structure of type $A\in {\cal M}$. 
In fact, there are only two such 
groups which are described in detail
in Section \ref{resum}. We show
\begin{prop}
There are two finite groups which admit
a mixed ramification structure of type $A\in {\cal M}$. They both
have order $256$ and admit a structure of type $[4,4,4]_{16}\in{\cal M}_3$.
\end{prop}

The proof  relies again heavily on the use of the MAGMA-library containing
all groups of low order. 
In order to avoid an excess of computations
we first consider each type $A=[m_1,\ldots,m_r]\in {\cal M}_r$ 
($r\in\BN$) seperately going through
the finite list of Proposition \ref{propmtup}, 
trying in a first round to exclude as
many cases as possible by some criteria, which are computationally 
cheap to verify. 

If these are satisfied and if we have we have access to the groups of order
$\beta(A)^2/2$ through a MAGMA-library we check for each of these groups $H$
\begin{itemize}
\item does $H$ admit a spherical system of generators of type $A$?
\item does $H$ admit a disjoint pair of spherical systems of 
generators of type $(A,A)$?
\end{itemize}
In fact, only very few groups $H$ survive the first test and
 for them the second
criterion, though computationally expensive, can be carried out. We are left
with a small list of groups $H$ admitting a disjoint pair of spherical systems
 of generators of type $(A,A)$, i.e., an unmixed ramification structure. 
Fortunately such groups $H$
only appear when the order $\beta(A)^2$ is small enough to have access to
all groups $G$ of this order. We then go through all these groups $G$ and 
list their subgroups of index $2$ isomorphic to one of the groups $H$. We
then check whether the compatibility conditions for a mixed  
ramification structure of type $A$ 
(see Section \ref{groups}) could be satisfied.
   
If the order  $\beta(A)^2/2$ is too big to use a 
MAGMA-library we analyse the 
subgroups of low index in the polygonal group 
$\BT(m_1,\ldots,m_r)$ 
(which often happen to be isomorphic to  
polygonal groups). We always can show that $H$ would then have
a subgroup of low index which is a quotient of another polygonal 
group. Using this descent procedure, sometimes repeatedly, 
always brought us into a
region of orders accessible to MAGMA-libraries.

\begin{rem}
Proposition $(4.4)$ of \cite{BCG} contains a misprint (the order of $G$
was mistakenly confused with the order of $H$). The correct
statement is: no group of order $< 256$ admits a mixed Beauville structure, i.e., a
mixed ramification structure of length $3$.
\end{rem}
%%%%%%%%%%%%%%%%%%%%%%%%%%%%%%%%%%%%%%%%%%%%%%%%%%%%%%%%
\subsection{$A\in{\cal M}_3$, $\beta(A)\ne 16$}
%%%%%%%%%%%%%%%%%%%%%%%%%%%%%%%%%%%%%%%%%%%%%%%%%%%%%%%%

In this section we treat the cases $A\in{\cal M}_3$, $\beta(A)\ne 16$. In each
case we indicate a sequence of MAGMA-computations showing that there is no 
finite group $G$ admitting
a mixed ramification structure of such type.

%%%%%%%%%%%%%%%%%%%%%%%%%%%%%%%%%%%%%%%%%%%%
\subsubsection{\bf $A\in{\cal M}_3$,\ $\beta(A)=8,\, 10,\, 12,\, 14,\, 18,
\, 20,\, 30,\, 36,\, 42,\, 60$}
%%%%%%%%%%%%%%%%%%%%%%%%%%%%%%%%%%%%%%%%%%%%

In these cases the order of $H$ is either small or divisible
only by a power of $2$ with exponent $\leq 3$ and
all relevant groups $H$ can be quickly inspected by MAGMA. 
We find
\begin{comp}
Let $k$ be one of the numbers $8,\, 10,$ $12,\, 14,$ 
$18,\, 20,$ $30,\, 36,$ $42,\, 60$
and $H$ a group of order $k^2/2$. Then
$H$ does not have a disjoint pair of spherical systems of 
generators of type $(A,A)$
with $\beta(A)=k$. 
\end{comp}
For $k=10,\, 14,\, 60$ there is even no group of order  $k^2/2$ having a 
system of generators of type $A$
with $\beta(A)=k$.

%%%%%%%%%%%%%%%%%%%%%%%%%%%%%%%%%%%%%%%%%%%%%%
\subsubsection{\bf $A\in{\cal M}_3$,\ $\beta(A)=24$}
%%%%%%%%%%%%%%%%%%%%%%%%%%%%%%%%%%%%%%%%%%%%%%

The group $G$ has order $|G|=576= 32 \cdot 27$, hence is
solvable. Since the order of the subgroup $H$ is still low we may quickly infer
from MAGMA
\begin{comp}
No group $H$ of order $288$
has a disjoint pair of spherical systems of generators of type $(A,A)$
with $A\in{\cal M}_3$ and $\beta(A)=24$. 
\end{comp}

%%%%%%%%%%%%%%%%%%%%%%%%%%%%%%%%
\subsubsection{\bf $A\in{\cal M}_3$,\ $\beta(A)=32$}
%%%%%%%%%%%%%%%%%%%%%%%%%%%%%%%%

Here $A=[ 2 , 4 ,  8]_{ 32 }$.
The group $G$ has order $|G|=1024$, it is a $2$-group.
Its subgroup $H$ has order $512$ and 
has a spherical system of generators of type $[ 2 , 4 ,  8]$. A  computation (of about $6$
hours) 
using the MAGMA-library of groups of order $512$ reveals
\begin{comp}
There are $10494213$ groups of order $512$. Eight of them
($\SmG(512,v)$ for $v=409$, $1818$, $1822$, $1832$, $1838$, $1854$, $1862$,
$2023$)   
have a system of generators of type $[ 2 , 4 ,  8]$.
\end{comp}
\medskip
The analysis of so many groups is made feasible by first selecting those
groups of order $512$ whose abelianisation is a quotient of the 
abelianisation of $\BT(2,4,8)$ which is isomorphic to $\Z_2\times\Z_4$.

Now a quick computation shows
\begin{comp}
None of the groups of order $512$ admitting  
a spherical system of generators of type $[ 2 , 4 ,  8]$ has a 
disjoint pair of spherical systems of generators of type 
$([ 2 , 4 ,  8],\,[ 2 , 4 ,  8])$. 
\end{comp}

%%%%%%%%%%%%%%%%%%%%%%%%%%%%%%%%%%%%%
\subsubsection{\bf $A\in{\cal M}_3$,\ $\beta(A)=40$}
%%%%%%%%%%%%%%%%%%%%%%%%%%%%%%%%%%%%%

Here $A=[ 2 , 5 ,  5]_{40}$. 
The group $G$ has order $|G|=1600$, hence it is
solvable. 
Its subgroup $H$ has order $800$ and has a spherical system of generators 
of type  $[2,5,5]$.
We have 
\begin{equation}
\BT(2,5,5)^{\rm ab}=\Z_5.
\end{equation}
Therefore the abelianisation of $H$ is $\Z_5$.
\begin{comp} 
There are $1211$ groups of order $800$. 
None of them has abelianisation $\Z_5$.
\end{comp}

%%%%%%%%%%%%%%%%%%%%%%%%%%%%%%%%%%%%%%%
\subsubsection{\bf $A\in{\cal M}_3$,\ $\beta(A)=48$}
%%%%%%%%%%%%%%%%%%%%%%%%%%%%%%%%%%%%%%%%

Here we have
$A=[2,3,12]_{48},\, [ 2 , 4 ,  6]_{ 48},\,  [ 3 , 3 ,  4]_{ 48}$.
Going through the MAGMA-library of groups of order $1152=24\cdot 48$
we find 
\begin{comp} 
There are $157877$ groups of order $1152$. 
\begin{itemize}
\item None of them has a system of generators of type $[ 3,3 ,4]$, 
\item one of them
($\SmG(1152,155454)$) 
has a system of generators of type $[ 2 , 3 ,12]$, 
\item one of them ($\SmG(1152,157849)$)
has a system of generators of type $[ 2 , 4 ,6]$.
\end{itemize}
\end{comp}
Treating the two groups $\SmG(1152,155454)$ and
$\SmG(1152,157849)$ which are remaining, is computationally cheap. We have 
\begin{comp}  
1) If $[a,b,c]$, $[a',b',c']$
are two systems of generators the group $H:=\SmG(1152,155454)$ of type 
$[ 2 , 3 ,  12]$ then $a,\, a'$ are conjugate in $H$.

\noindent
2) If $[a,b,c]$, $[a',b',c']$
are systems of generators of $H:=\SmG(1152,157849)$ of type 
$[ 2 , 4 ,6]$ then $a,\, a'$ are conjugate in $H$.
\end{comp}
This shows that there are no groups $H$ of order $1152$ with a 
disjoint pair of spherical systems
of generators of one of the types above. 

%%%%%%%%%%%%%%%%%%%%%%%%%%%%%%%%%%%%%%%
\subsubsection{\bf $A\in{\cal M}_3$,\ $\beta(A)=72$}
%%%%%%%%%%%%%%%%%%%%%%%%%%%%%%%%%%%%%%%

Here $A=[2,3,9]_{72}$. The group $G$ has order $|G|=5184$, hence it is
solvable. Its index $2$
subgroup $H$ has order $2592$ and system of generators of type  $[2,3,9]$.
The order of $H$ is too big to apply the computational arguments used 
before. To exclude this case we argue as follows.

We have 
\begin{equation}\label{we5}
\BT(2,3,9)^{\rm ab}=\Z_3,\qquad \BT(2,3,9)' \cong \BT(2,2,2,3).
\end{equation}
The commutator subgroup $H'$ of $H$ has order $864$ and we conclude from 
(\ref{we5}) 
that $H'$ is a quotient of $\BT(2,2,2,3)$. Looking through the
relevant MAGMA-library we find
\begin{comp}\label{compufa72}
Of the $4725$ groups of order $864$ only the two groups 
$H_1:=\SmG(864,2225)$ and
$H_2:=\SmG(864,4175)$ are quotients of $\BT(2,2,2,3)$.
\end{comp}
We are left with the question whether $H_1,\, H_2$ are the commutator 
subgroup of a group $H$ of order $2592$ which has a disjoint pair of spherical 
systems of generators of type $(A,A)$. We need
\begin{lemma}\label{lemgru72}   
Let $\tilde H$ be a finite group with $\tilde H^{\rm ab}=\Z_3$
which has a disjoint pair of spherical 
systems of generators
$([h_{(1,1)},h_{(1,2)},h_{(1,3)}],\,[h_{(2,1)},h_{(2,2)},h_{(2,3)}])$  with $[2,3,9]=$
$[{\rm ord}(h_{(1,1)}),{\rm ord}(h_{(1,2)}), {\rm ord}(h_{(1,3)})]$
$=[{\rm ord}(h_{(2,1)}),{\rm ord}(h_{(2,2)}), {\rm ord}(h_{(2,3)})]$. Then
\begin{itemize}
\item[{\rm (i)}]
$[h_{(i,1)},h_{(i,2)}h_{(1,1)}h_{(i,2)}^{-1},h_{(i,2)}^2h_{(i,1)}h_{(i,2)}^{-2}]$
($i=1,\, 2$) are generating tuples of elements of order $2$ for $\tilde H'$ such that 
$z_i:=h_{(i,1)}\cdot h_{(i,2)}h_{(1,1)}h_{(i,2)}^{-1}\cdot 
h_{(i,2)}^2h_{(i,1)}h_{(i,2)}^{-2}$ has order $3$.
 
\item[{\rm (ii)}] $\tilde \Sigma_1\cap\tilde \Sigma_2=\emptyset$ where for
$i=1,\, 2$
$$\tilde\Sigma_i:=\bigcup_{h\in \tilde H'} h\{\,h_{(i,1)},\,    
z_i,\, z_i^2\,\}h^{-1}$$
\end{itemize}
\end{lemma}
We skip the straightforward proof.

We finish this case by
\begin{comp} The groups $H_1$, $H_2$ in Computational Fact \ref{compufa72}
do not have generating $3$-tuples with
properties ${\rm (i)},\, {\rm (ii)}$ of Lemma \ref{lemgru72}.
\end{comp}

%%%%%%%%%%%%%%%%%%%%%%%%%%%%%%%%%%%%%%%%%%%%%%%%
\subsubsection{\bf $A\in{\cal M}_3$,\ $\beta(A)=80$}
%%%%%%%%%%%%%%%%%%%%%%%%%%%%%%%%%%%%%%%%%%%%%%%%

Here $A=[2,4,5]_{80}$. 
The group $G$ has order $|G|=6400$, hence is
solvable. Its index $2$ subgroup $H$ has a system of 
generators of type  $[2,4,5]$.
We have 
\begin{equation}
\BT(2,4,5)^{\rm ab}=\Z_2,\qquad \BT(2,4,5)' \cong \BT(2,5,5).
\end{equation}
This implies that there must be a subgroup $H_1$ of index
$2$ in $H$ which is a quotient of $\BT(2,5,5)$.
$H_1$  has order $1600$ and is solvable, hence it has abelianisation 
$\Z_5$. An inspection of the relevant MAGMA-library shows
\begin{comp}
There are $10281$ groups of order $1600$. 
None of them has abelianisation $\Z_5$. 
\end{comp}

%%%%%%%%%%%%%%%%%%%%%%%%%%%%%%%%%%%%%%%%
\subsubsection{\bf $A\in{\cal M}_3$,\ $\beta(A)=96$}
%%%%%%%%%%%%%%%%%%%%%%%%%%%%%%%%%%%%%%%%

Here $A=[2,3,8]_{96}$ and
$G$ has order $|G|=96^2=9216$, hence is
solvable. Its subgroup $H$ has a system of generators of type  $[2,3,8]$.

We have 
\begin{equation}
\BT(2,3,8)^{\rm ab}=\Z_2,\qquad \BT(2,3,8)' \cong \BT(3,3,4).
\end{equation}
This implies that there must be a subgroup $H_1$ of index
$2$ in $H$ which is a quotient of $\BT(3,3,4)$. We further have
\begin{equation}
\BT(3,3,4)^{\rm ab}=\Z_3, \qquad  \BT(3,3,4)'\cong \BT(4,4,4).
\end{equation}
This in turn implies ($H$ is solvable) that there is a subgroup $H_2$ of index
$3$ in $H_1$ which is a quotient of $\BT(4,4,4)$. Note that $|H_2|=768$. 
\begin{comp}
There are $1090235$ groups of order $768$. 
None of them is a quotient of $\BT(4,4,4)$. 
\end{comp}
This fact can be derived by inspection of the relevant MAGMA-library. 
The analysis of such a hugge number of groups is made possible by first selecting
those groups of order $768$ which have abelianisation which is a quotient of the 
abelianisation of $\BT(4,4,4)$ which is $\Z_4\times\Z_4$. There are $1651$ such
groups. For them it is quickly checked whether they are a quotient of
$\BT(4,4,4)$.

%%%%%%%%%%%%%%%%%%%%%%%%%%%%%%%%%%%%%%%%%%%%%%
\subsubsection{\bf $A\in{\cal M}_3$,\ $\beta(A)=168$}
%%%%%%%%%%%%%%%%%%%%%%%%%%%%%%%%%%%%%%%%%%%%%%

Here $A=[2,3,7]_{168}$. 
The group $G$ has order $|G|=28224$, its subgroup $H$ has order
$14112$ and has a system of generators of type $[2,3,7]$.

The group $H$ is perfect, hence has a non abelian simple group as quotient. 
Note that the only non abelian simple groups with order dividing $14112$
are $\PSL(2,\BF_7)$ and $\PSL(2,\BF_8)$. 
Let $K$ be the kernel of the
quotient homomorphism. This group has order $28$ or $84$ and 
it has a normal, hence characteristic
Sylow-$7$-subgroup $S$. 

Therefore $S$ has to be normal in $G$ and $G/S$ has
order $2016$ and is perfect. There are no such groups.

%%%%%%%%%%%%%%%%%%%%%%%%%%%%%%%%%%%%%%%%%%%%%%%%%%%%%%%%
\subsection{$A\in{\cal M}_3$, $\beta(A)=16$}\label{resum}
%%%%%%%%%%%%%%%%%%%%%%%%%%%%%%%%%%%%%%%%%%%%%%%%%%%%%%%%

In Section \ref{tup} we have shown that the possible tuples are
$[ 2 , 6 ,  12]_{ 16 }$, 
$[ 3 , 3 ,  12]_{ 16 }$, $[ 3 , 4 ,  6]_{ 16 }$, $[ 2 , 8 ,  8 ]_{ 16 }$
$[ 4 , 4 ,  4]_{ 16 }$. The first three are not possible since a group of
order $128$ cannot contain elements of order $3$. This leaves 
$[ 2 , 8 ,  8 ]_{ 16 }$ and $[ 4 , 4 ,  4]_{ 16 }$.

\subsubsection{\bf $A=[ 2 , 8 ,  8 ]_{ 16 }$}

Going through the list of groups of order $128$ we find
\begin{comp}
There are $2328$ groups of order $128$. 
Only $7$ of them have a spherical system of
generators of type $[ 2 , 8 ,  8 ]$. These are
$\SmG(128,v)$ for $v=2,\, 48,\, 50,\, 77,\, 135,\, 137,\, 142$.
\end{comp}
Analysing the $7$ remaining groups we find
\begin{comp}
None of the groups $\SmG(128,v)$ ($v=2$,\, $48,$\, $50,$\, $77,$\, 
$135,$\, $137,$\, $142$)
has a disjoint pair of spherical systems of generators of type $(A,A)$.
\end{comp}

%%%%%%%%%%%%%%%%%%%%%%%%%%%%%%%%%%%%%%%%%%%%%%%%%%
\subsubsection{\bf $A=[ 4 , 4 ,  4 ]_{ 16 }$}\label{suse256}
%%%%%%%%%%%%%%%%%%%%%%%%%%%%%%%%%%%%%%%%%%%%%%%%%%

Going through the list of groups of order $128$ we find
\begin{comp}
There are $2328$ groups of order $128$. Only $4$ of them have a system of
generators of type $[ 4 , 4 ,  4 ]$. These are
$\SmG(128,v)$ for $v=36,\, 125,\, 141,\, 144$. Only 
$\SmG(128,36)$ has  a disjoint pair of systems of generators of type $(A,A)$.
\end{comp}
This leaves us with the possibility that $H\cong \SmG(128,36)$. Going through
the groups of order $256$ we find
\begin{comp}
Of the $56092$ groups of order $256$ only $29$ contain a subgroup of index $2$
isomorphic to $\SmG(128,36)$. They are the groups $\SmG(128,v)$ for $v=$ 
$ 382 $, $ 414 $, $ 1087 $, $ 1088 $, $ 1089 $, $ 1090 $, $ 1734 $, $ 1735 $,
$ 1736 $, $ 1737 $, $ 1738 $, $ 2483 $, $ 2484 $, $ 2485 $, $ 2486 $,
$ 2487 $, $ 2488 $, $ 2489 $, $ 2490 $, $ 3324 $, $ 3325 $, $ 3326 $, $ 3327 $,
$ 3378 $, $ 3379 $, $ 3380 $, $ 3381 $, $ 3678 $, $ 3679 $.
\end{comp}
Analysing the $29$ remaining groups we easily find
\begin{prop}
Of the groups of order $256$ exactly
\begin{equation}
{\rm G}(256,1):=\SmG(256,3678),\quad {\rm G}(256,2):=\SmG(256,3679)
\end{equation}
admit a mixed ramification structure of type $[ 4 , 4 ,  4 ]_{ 16 }$.
\end{prop}

Presentations for the groups
${\rm G}(256,1)$ and ${\rm G}(256,2)$ are
\begin{equation}\label{gr2561}
{\rm G}(256,1)=\left\langle\, g_1,\ldots,g_8\Mid  
\begin{matrix}
g_1^2 = g_4g_5g_6, & g_2^2 = g_4  g_5, &  g_3^2 = g_4, \cr
g_2^{g_1} = g_2  g_4, & g_3^{g_1} = g_3  g_5, & g_3^{g_2} = g_3  g_6,\cr
g_4^{g_1} = g_4  g_7, & g_4^{g_2} = g_4  g_8, & g_5^{g_1} = g_5g_7 g_8,\cr
g_5^{g_2} = g_5  g_8, & g_5^{g_3} = g_5  g_7, & g_6^{g_1} = g_6  g_8,\cr
g_6^{g_2} = g_6  g_7, & g_6^{g_3} = g_6  g_8 &
\end{matrix}
\right\rangle,
\end{equation}
\begin{equation}\label{gr2562}
{\rm G}(256,2)=\left\langle\, g_1,\ldots,g_8 \Mid  
\begin{matrix}
g_1^2 = g_4g_5g_6g_7, & g_2^2 = g_4  g_5, &  g_3^2 = g_4, \cr
g_2^{g_1} = g_2  g_4, & g_3^{g_1} = g_3  g_5, & g_3^{g_2} = g_3  g_6,\cr
g_4^{g_1} = g_4  g_7, & g_4^{g_2} = g_4  g_8, & g_5^{g_1} = g_5g_7 g_8,\cr
g_5^{g_2} = g_5  g_8, & g_5^{g_3} = g_5  g_7, & g_6^{g_1} = g_6  g_8,\cr
g_6^{g_2} = g_6  g_7, & g_6^{g_3} = g_6  g_8 &
\end{matrix}
\right\rangle.
\end{equation}
The conventions for these so called PC-presentations are explained in Section
\ref{suse44}.  
       
%%%%%%%%%%%%%%%%%%%%%%%%%%%%%%%%%%%%%%%%%%%%%%%%%%%%%%%%
\subsection{$A\in{\cal M}_4,\ {\cal M}_5,\ {\cal M}_6,\ {\cal M}_8$}
%%%%%%%%%%%%%%%%%%%%%%%%%%%%%%%%%%%%%%%%%%%%%%%%%%%%%%%

In this section we treat the cases $A\in{\cal M}_4,\, 
{\cal M}_5,\, {\cal M}_6,\, {\cal M}_8$. We show
\begin{prop} There is no finite group $G$ admitting a mixed ramification
structure of type $A\in{\cal M}_4,\, 
{\cal M}_5,\, {\cal M}_6,\, {\cal M}_8$.
\end{prop}

%%%%%%%%%%%%%%%%%%%%%%%%%%%%%%%%%%%%%%%%
\subsubsection{$A\in{\cal M}_4$}
%%%%%%%%%%%%%%%%%%%%%%%%%%%%%%%%%%%%%%%%

The order of the group $H$ is at most $288$ and all relevant groups can be
checked for generating systems. 

\begin{center}
\begin{tabular}{|c|c|}
\hline
A  &   \\
\hline
$[  2 ,  2 ,  2 ,  3]_{24}$ & No disjoint generating systems\\
$[  2 ,  2 ,  2 ,  4 ]_{16 }$ & No disjoint generating systems\\ 
$[  2 ,  2 ,  2 ,  6 ]_{12 }$ & No disjoint generating systems\\ 
$[  2 ,  2 ,  3 ,  3 ]_{  12 }$ & No disjoint generating systems\\ 
$[  2 ,  2 ,  2 ,  10 ]_{10}$ & No generating systems\\
$[  2 ,  2 ,  4 ,  4 ]_{8}$ & $\SmG(32,22)$ admits a disjoint generating systems\\ 
$[  2 ,  2 ,  6 ,  6 ]_{6 }$ & No disjoint generating systems\\
$[  2 ,  3 ,  3 ,  6 ]_{6 }$ & No disjoint generating systems\\
$[  3 ,  3 ,  3 ,  6 ]_{6 }$ & No generating systems\\
$[  4 ,  4 ,  4 ,  4 ]_{4 }$ & No disjoint generating systems\\
\hline
\end{tabular}
\end{center}

\begin{comp}
Of the $267$ groups of order $64$ only $32$ contain a subgroup of index $2$
which is isomorphic to $\SmG(32,22)$. None of them has a mixed ramification
structure of type $[  2 ,  2 ,  4 ,  4 ]_{8}$.
\end{comp}

\subsubsection{$A\in {\cal M}_5,\ {\cal M}_6,\ {\cal M}_8$}

In these cases the relevant group orders are so small that all groups $G$ can
easily be inspected. There is none with a mixed ramification structure.

%%%%%%%%%%%%%%%%%%%%%%%%%%%%%%%%%%%%%%%%%%%%%%%%
\section{Moduli spaces}\label{modu}
%%%%%%%%%%%%%%%%%%%%%%%%%%%%%%%%%%%%%%%%%%%%%%%%

In this section we will describe completely the moduli spaces of the surfaces
isogenous to a product with $p_g = q = 0$. More precisely, let
$\mathfrak{M}_{(1,8)}$ be the moduli space of minimal smooth complex
projective surfaces with $\chi (S) = 1$ and $K_S^ 2 = 8$. As usual $K_S$
denotes the canonical divisor of $S$ and $\chi(S) = 1 + p_g(S) -q(S)$ is the
holomorphic Euler-Poincare' characteristic of $S$. It is nowadays wellknown
(cf. \cite{gieseker}) that $\mathfrak{M}_{(a,b)}$ is quasiprojective for all
$a,\, b\in\BN$.
Obviously, our surfaces are contained in the 
moduli space $\mathfrak{M}_{(1,8)}$
and we will
describe their locus there.

Let $G$ be a finite group and fix an unmixed ramification type 
$(A,B)\in \BN^r\times \BN^s$. We
denote by $\mathfrak{M}_{(G;A,B)}$ the subset of $\mathfrak{M}_{(1,8)}$ defined
by isomorphism classes of surfaces isogenous to a product admitting
a  ramification type
$(A,B)$ (or $(B,A)$).

We observe
\begin{rem} 1) The set
$\mathfrak{M}_{(G;A,B)}\subset \mathfrak{M}_{(1,8)}$ 
consists of a finite number of connected components
 of the same dimension, 
which are irreducible in the Zariski topology.

2) It is clear from Section \ref{basi} that the dimension 
$d(G;A,B)$ of any component in
$\mathfrak{M}_{(G;A,B)}$ is precisely $\ell(A)-3+\ell(B)-3$
since we take $\ell(A)$-points in $\PP^1$ modulo projective equivalence,
and likewise $\ell(B)$-points in $\PP^1$ modulo projective equivalence.
\end{rem}

In order to calculate the number of components $n(G;A,B)$ of
$\mathfrak{M}_{(G;A,B)}$ we use the following
\begin{prop}\label{thu1}
Let $S,\, S'$ be a surfaces isogenous to a product, 
of unmixed type and with $q(S)=q(S')=0$. Then $S$, $S'$ 
are in the same irreducible component if and only if $G(S)\cong G(S')$, 
$(A_1(S),A_2(S))=(A_1(S'),A_2(S'))$ and
${\cal T}(S)$ and ${\cal T}(S')$ are in the same orbit of 
 ${\bf B}_r\times {\bf B}_s\times {\rm Aut}(G)$ where $r=\ell(T_1)$, 
$s=\ell(T_2)$.
\end{prop}
For a proof we refer to \cite{BaCa}.

By computer calculation we obtain the following table 
of the possible unmixed ramification structures on finite groups 
of type $(A,B)$ with $A,\, B\in {\cal N}$ leading to surfaces with $K^2=8$.

\begin{theo}
If $S \neq \PP^1 \times \PP^1$ is a smooth projective surface isogenous to a
product of unmixed type with $p_g(S)=q(S)=0$ 
and
with minimal realisation $S\cong (C_1\times C_2) /G$ then 
$G$ is one of the
groups in the following table and the genera of the curves $C_1,\, C_2$ are as listed   
in the table. The numbers of components $N$ in $\mathfrak{M}_{(1,8)}$ and
their dimension is given in the remaining two columns. 
\begin{center}
\begin{tabular}{|c|c|c|c|c|c|}
\hline
$G$  & $|G|$ & $A$ & $B$ & $n(G;A,B)$ & $d(G;A,B)$ \\
\hline
${\mathfrak A}_5$ & $60$ & $[2,5,5]_{20}$ & $[3, 3, 3, 3]_{3}$  & 1 & 1\\
${\mathfrak A}_5$ & $60$ & $[5,5,5]_{ 5 }$ & $[ 2, 2, 2, 3]_{12}$ & 1& 1\\
${\mathfrak A}_5$ & $60$ & $[3,3,5]_{15}$ & $[2,2,2,2]_{4}$ & 1& 1\\
${\mathfrak S}_4\times \Z_2$ & $48$ & $[2,4,6]_{24}$ & $[2,2,2,2,2,2]_2$ & 1 &
3\\ G($32$) & $32$ & $[ 2,2,4,4]_4$ & $[2,2,2,4]_8$ & 1 & 2\\
$\Z_5^2$ & $25$ & $[5,5,5]_5$ & $[5,5,5]_5$ & 2 & 0 \\
${\mathfrak S}_4$ & $24$ & $[ 3, 4, 4]_{12}$ & $[2,2,2,2,2,2]_2$ & 1& 3\\
G($16$) & $16$ & $[ 2, 2, 4, 4]_4$ & $[2,2,4,4]_4$ & 1& 2\\
${\rm D}_4\times\Z_2$ & $16$ & $[2,2,2,4]_8$ & $[2,2,2,2,2,2]_2$ & 1& 4\\
$\Z_2^4$ & $16$ & $[2,2,2,2,2]_4$ & $[2,2,2,2,2]_4$ & 1 & 4 \\
$\Z_3^2$ & $9$ & $[ 3, 3, 3, 3]_3$ & $[3,3,3,3]_3$ & 1 & 2 \\ 
$\Z_2^3$ & $8$ & $[2,2,2,2,2]_4$ & $[2,2,2,2,2,2]_2$ & 1 & 5 \\ 
\hline
\end{tabular}
\end{center}
\end{theo}

\medskip
The case $G$ abelian was done in \cite{BaCa}. In \cite{PardDP} four of the non
abelian cases are constructed and for three of these the irreducibility of the
corresponding family is proven.

We now turn to the mixed case.

Let $G$ be a finite group and fix a mixed ramification type $A$. We
denote by $\mathfrak{M}_{(G;A)}$ the subset of $\mathfrak{M}_{(1,8)}$ 
given
by the isomorphism classes of surfaces isogenous to a product 
admitting a  mixed
ramification type $A$.

Also here we have
\begin{rem} The set 
$\mathfrak{M}_{(G;A)}\subset \mathfrak{M}_{(1,8)}$ consists of a 
finite number of connected components
 of the same dimension 
$d(G;A) = \ell(A) - 3$, which are irreducible in
the Zariski topology.
\end{rem}

\begin{prop}
Let $S,\, S'$ be surfaces isogenous to a product, 
of mixed type and with $q(S)=q(S')=0$. Then $S$, $S'$ 
are in the same irreducible component if and only if $G(S)\cong G(S')$ and
${\cal T}(S)$ and ${\cal T}(S')$ are in the same orbit of 
 ${\bf B}_r\times {\rm Aut}(G)$ where $r=\ell(T)$.
\end{prop}

Hence the number of components $n(G;A)$ of $\mathfrak{M}_{(G;A)}$ is
precisely the number of orbits of ${\bf B}_{\ell(A)} \times {\rm Aut}(G)$ 
on the set $\mathcal{B}(G;A)$.

We already know from Section \ref{classimi} that there are exactly the 
groups ${\rm G}(256,1)$, ${\rm G}(256,2)$ which have a 
mixed ramification structure
of type $A\in {\cal M}$. In fact, they both have such a structure 
of type $A=[4,4,4]_{16}$.
We shall now determine the numbers of orbits of  
${\bf B}_3\times {\rm Aut}({\rm G}(256,i)$ ($i=1,\, 2$) on the set of
ramification structures.

Let us begin with $G={\rm G}(256,1)$. We have
\begin{prop}\label{prop2561}
\begin{itemize}
\item[{\rm (i)}] The automorphism group 
of ${\rm G}(256,1)$ has $12288$ elements, it acts
with $3$ orbits on the set of subgroups of index $2$ in ${\rm G}(256,1)$.
\item[{\rm (ii)}] Representatives for the $3$ orbits are 
$H_1:=\langle \, g_1,\, g_3 \,\rangle$ (which is a fixed point for 
${\rm Aut}({\rm G}(256,1))$), 
$H_2:=\langle \, g_1,\, g_2 \,\rangle$ (which has an orbit of cardinality $3$,
$H_3:=\langle \, g_2,\, g_1g_3 \,\rangle$ 
(which has an orbit of cardinality $3$).
\item[{\rm (iii)}] The action of ${\bf B}_3\times {\rm Aut}({\rm G}(256,1)$ on 
$\mathcal{B}({\rm G}(256,1);[4,4,4]_{16})$ has $3$ orbits (corresponding to
the $3$ orbits of ${\rm Aut}({\rm G}(256,1))$ on 
the set of subgroups of index $2$ in ${\rm G}(256,1)$. 
\end{itemize} 
\end{prop}

\medskip
For the group $G={\rm G}(256,2)$ the picture is different, we find
\begin{prop}
\begin{itemize}
\item[{\rm (i)}] The automorphism group 
of ${\rm G}(256,2)$ has $86016$ elements, it acts
transitively on the set of subgroups of index $2$ in ${\rm G}(256,2)$.
\item[{\rm (ii)}] The action of ${\bf B}_3\times {\rm Aut}({\rm G}(256,2))$ on 
$\mathcal{B}({\rm G}(256,2);[4,4,4]_{16})$ is transitive.
\end{itemize}
\end{prop}

\medskip
The proof of the above two propositions is done by standard MAGMA routines.
s
Combining these results we find the following table 
of the possible mixed ramification structures on finite groups 
of type $A$ with $A\in {\cal M}$.

\begin{theo}
If $S$ is a smooth projective surface isogenous to a
product of mixed type with $p_g(S)=q(S)=0$ 
and
with minimal realisation $S\cong (C_1\times C_2) /G$, then 
$G$ is one of the
groups in the following table and the genera of the curves $C_1,\, C_2$ are as listed   
in the table. The numbers $N : = n(G;A)$ of components of $\mathfrak{M}_{(1,8)}$
and their dimension is given in the remaining two columns. 
\begin{center}
\begin{tabular}{|c|c|c|c|c|}
\hline
$G$  & $|G|$ & $A$ & $n(G;A)$ & $d(G;A)$ \\
\hline
$G={\rm G}(256,1)$ & $256$ & $[4,4,4]_{16}$ & 3 & 0\\
$G={\rm G}(256,2)$ & $256$ & $[4,4,4]_{16}$ & 1 & 0\\
\hline
\end{tabular}
\end{center}
\end{theo}
%%%%%%%%%%%%%%%%%%%%%%%%%%%%%%%%%%%%%%%%%%%%%%%%
\section{Concrete models}\label{concrete}
%%%%%%%%%%%%%%%%%%%%%%%%%%%%%%%%%%%%%%%%%%%%%%%%

In this section we want to give explicit descriptions of the groups
and spherical systems of generators occurring in the nonabelian case
(the abelian case is fully classified and described in \cite{BaCa}).

Some of these nonabelian examples were already described in \cite{BaCa}, 
but we thought it would be worthwhile to give a complete list.

\subsection{ $G = \mathfrak A_5$ }

The unmixed ramification structure of type 
$$ ([ 3,3,3,3 ] , [2,5,5]) $$
is given by the following elements of $ \mathfrak A_5$: 
$$ ([(1,2,3),(3,4,5), (4,3,2), (2,1,5)],[(2,4)(3,5), (2,1,3,4,5), (1,2,3,4,5)]).$$

The unmixed ramification structure of type 
$$ ([ 5,5,5 ] , [2,2,2,3]) $$
is given by the following elements of $ \mathfrak A_5$: 
$$ ([(1,2,5,3,4), (1,2,4,5,3),(1,2,3,4,5) ],  [(1,2)(3,4), (2,4)(3,5), (1,4)(3,5),
(2,3,4)]).$$

The unmixed ramification structure of type 
$$ ([ 2,2,2,2,2 ] , [3,3,5]) $$
is given by the following elements of $ \mathfrak A_5$: 
$$ ([(1,2)(3,4), (1,3)(2,4),  (1,4) (2,3), (1,4)(2,5), (1,4)(2,5)],  [(1,2,3), (3,4,5) ,
(5,4,3,2,1)]).$$

\subsection{ $G = \mathfrak D_4 \times \Z_2$ }

We write as customary $\mathfrak D_4$ as the group generated by elements $x, y$
satisfying the relations $ x^4 = y^2 = 1, yxy = x^{-1}.$

Then there is exactly one class of unmixed ramification structures, of type 
$$ ([ 2,2,2,2,2 ,2] , [2,2,2,4]) $$
 given by the following elements of $ \mathfrak D_4 \times \Z_2$ 

$$ ([(y,0), (yx,1), (y x^2, 0), (yx,1), (x^2,1), (x^2,1)], [(1,1), (y,1), (xy,0), (x,0)]). $$

\subsection{ $G = \mathfrak S_4 $ }

There is exactly one class of unmixed ramification structures, of type 
$$ ([ 2,2,2,2,2 ,2] , [3,4,4]) $$
 given by the following elements of $  \mathfrak S_4 $:

$$([(1,2), (1,2) , (2,3), (2,3), (3,4),  (3,4)],  [(1,2,3), (1,2,3,4), (1,2,4,3)]). $$

Note that this generating system is contained in the ArXiv version of \cite{BaCa},
it was not possible for technical reasons to correct the
printed version in time.

\subsection{ $G = \mathfrak S_4  \times \Z_2$ }

There is exactly one class of unmixed ramification structures, of type 
$$ ([ 2,4,6] , [2,2,2,2,2,2]) $$
 given by the following elements of $  \mathfrak S_4 \times \Z_2$:

$$([[(1,2),0], [ (1,2,3,4),1], [ (4,3,2),1]] ,$$
$$ [[(1,2) (3,4), 1], [ (1,2),1], [(3,4),1], 
[(2,3) (1,4), 1], [ (2,3),1], [(1,4),1]]). $$

\subsection{ $G = G(16)$.}

We use here the following realization of  $G : = {\bf G(16)}$
as a semidirect product 

$$ (\Z_4 \times \Z_2)\rfish_{\Phi} \Z_2$$ 
generated by $x,y,z$,  with centre $C \cong \Z_2 \times \Z_2$
generated by $x^2, y$, and such that 
$$ zxz = xy. $$

There is exactly one class of unmixed ramification structures, of type 
$$ ([ 2,2,4,4] , [2,2,4,4]) $$
 given by the following elements of ${\bf G(16)}$:

$$ ([z , z, x, x^{-1}], [z x^2 y,z x^2 y, xyz ,  (xyz)^{-1}]) .$$

\subsection{ $G = G(32),  G(256,1), G(256,2).$ }

We construct now concrete models for the finite groups like
$G(256,1)$ which make hand computations simple. We start off by giving a
general construction principle for metabelian groups. A group is called
metabelian if it contains an abelian normal subgroup with abelian quotient.   

Let now $N,\, Q$ be two abelian groups written additively. Let 
\begin{equation}
\Phi: Q\to {\rm Aut}(N), \qquad \Phi: q\mapsto \Phi_q \ \ (q\in Q)
\end{equation}
a homomorphism from $Q$ to the automorphism group of $N$. Further let
\begin{equation}
\Theta : Q\times Q\to N
\end{equation}
be a bilinear map. We define a multiplication on the set $N\times Q$ by 
setting
\begin{equation}
(n_1,q_1)\cdot (n_2,q_2):=
\left(n_1+\Phi_{q_1}(n_2)+\Theta(q_1,q_2),q_1+q_2\right)
\end{equation}
for $n_1,\, n_2\in N$ and $q_1,\, q_2\in Q$. We obtain a group structure iff 
$$\Phi_{q_1}(\Theta(q_2,q_3))=\Theta(q_2,q_3)$$
holds for all $q_1,\, q_2,\, q_3$. The resulting group is denoted by
\begin{equation}
N\rfish_{\Phi,\Theta} Q.
\end{equation}
There is the obvious exact sequence 
$$\langle 0\rangle\to N\to N\rfish_{\Phi,\Theta} Q\to Q\to \langle 0\rangle$$
hence $ N\rfish_{\Phi,\Theta} Q$ is metabelian. Conversely, every metabelian
group arises in this way. If $\Theta: Q\times Q\to N$ is the zero map then 
$N\rfish_{\Phi,\Theta} Q$ is a semidirect product of $N$ and $Q$ which we
denote by  $N\rfish_{\Phi} Q$.

We shall now describe the models for the remaining finite groups from Sections
\ref{classiumi} and \ref{classimi}.
 
\bigskip

\centerline{\bf G(32):}
\medskip

This group has nilpotency class $2$ and is a semidirect product of $\Z_2^4$ by 
$\Z_2$. The homomorphism 
$\Phi :\Z_2\to {\rm Aut}(\Z_2^4)={\rm GL}(4,\BF_2)$ can be
given by the single matrix (unipotent of order $2$). We set  
$$
\Phi_1:=\left(
\begin{matrix}
1 & 0 & 0 & 0\cr
0 & 1 & 0 & 0\cr
1 & 0 & 1 & 0\cr
0 & 1 & 0 & 1
\end{matrix}\right)
$$
From the presentation (\ref{gr32}) it can be seen that the 
resulting group $\Z_2^4\rfish_{\Phi} \Z_2$ is isomorphic to 
${\rm G}(32)$.

An unmixed ramification structure ${\cal T}=(T_1,T_2)$ of type
$([2,2,2,4]_8,[2,2,4,4]_4)$ 
on 
$\Z_2^4\rfish_{\Phi} \Z_2$ 
is given by 
$$T_1=[((0,0,1,1)^t,1),\ ((1,1,1,1)^t,0),\ 
((1,0,1,1)^t,0),\ ((0,1,1,1)^t,1)],$$
$$T_2=[((1,1,1,0)^t,0),\ ((1,0,0,0)^t,0),\ 
((1,1,1,0)^t,1),\ ((1,0,1,0)^t,1)].$$

\bigskip
\centerline{\bf G(256,1):}
\medskip

This group has nilpotency class $3$. But fortunately for us every group with
this property is metabelian. The group ${\rm G}(256,1)$ is 
of the form $\Z_2^5\rfish_{\Phi,\Theta}\Z_2^3$. We shall first describe
the maps $\Phi$ and $\Theta$.

Let $e_1,\, e_2,\, e_3$ be the standard basis of $\Z_2^3$. The homomorphism 
$\Phi :\Z_2^3\to {\rm Aut}(\Z_2^5)={\rm GL}(5,\BF_2)$ can be
given by its values on $e_1,\, e_2,\, e_3$. We set
$$
\Phi_{e_1}:=\left(
\begin{matrix}
1 & 0 & 0 & 0 & 0\cr
0 & 1 & 0 & 0 & 0\cr
0 & 0 & 1 & 0 & 0\cr
1 & 1 & 0 & 1 & 0\cr
0 & 1 & 1 & 0 & 1
\end{matrix}\right),
\quad
\Phi_{e_2}:=\left(
\begin{matrix}
1 & 0 & 0 & 0 & 0\cr
0 & 1 & 0 & 0 & 0\cr
0 & 0 & 1 & 0 & 0\cr
0 & 0 & 1 & 1 & 0\cr
1 & 1 & 0 & 0 & 1
\end{matrix}\right),
\quad
\Phi_{e_3}:=\left(
\begin{matrix}
1 & 0 & 0 & 0 & 0\cr
0 & 1 & 0 & 0 & 0\cr
0 & 0 & 1 & 0 & 0\cr
0 & 1 & 0 & 1 & 0\cr
0 & 0 & 1 & 0 & 1
\end{matrix}\right).
$$
To give the bilinear map $\Theta$ we set
$$\Theta(e_1,e_1):=(1,1,1,0,0)^t,\ \Theta(e_2,e_2):=(1,1,0,0,0)^t,\
\Theta(e_3,e_3):=(1,0,0,0,0)^t,$$
$$\Theta(e_2,e_1):=(1,0,0,1,1)^t,\
\Theta(e_3,e_1):=(0,1,0,0,1)^t,\
\Theta(e_3,e_2):=(0,0,1,1,1)^t$$
with the convention that the $\Theta(e_i,e_j)$ which are not mentioned are
equal to $0$.
From the presentation (\ref{gr2561}) it can be seen that the 
resulting group $\Z_2^5\rfish_{\Phi,\Theta} \Z_2^3$ is isomorphic to 
${\rm G}(256,1)$.

Here are three mixed ramification structures of type
$[4,4,4]_{16}$ 
on 
$\Z_2^5\rfish_{\Phi,\theta} \Z_2^3$:  
$$T_1:=[((0,0,1,0,1)^t,e_3),((1,1,0,0,0)^t,e_1),((1,0,0,1,0)^t,e_1+e3)],$$
$$T_2:=[((0,0,0,0,1)^t,e_2),((1,0,0,1,0)^t,e_1+e_2),((0,0,1,0,1)^t,e_1)],$$
$$T_3:=[((0,0,0,1,0)^t,e_2),((1,1,0,1,0)^t,e_1+e_2+e_3),
((1,0,1,1,0)^t,e_1+e_3)].$$
This is to say the three coordinates of $T_1,\, T_2,\, T_3$ generate a
subgroup $H$ of index $2$ in $G$ and the compatibility conditions of
Definition \ref{defimi} are satisfied. Moreover $T_1,\, T_2,\, T_3$ represent
the three orbits appearing in Proposition \ref{prop2561}. 

\bigskip
\centerline{\bf G(256,2):}
\medskip

This group has nilpotency class $3$ and is 
of the form $\Z_2^5\rfish_{\Phi,\Theta}\Z_2^3$. We shall describe
the maps $\Phi$ and $\Theta$.

Let $e_1,\, e_2,\, e_3$ be the standard basis of $\Z_2^3$. The homomorphism 
$\Phi :\Z_2^3\to {\rm Aut}(\Z_2^5)={\rm GL}(5,\BF_2)$ can be
given by its values on $e_1,\, e_2,\, e_3$. We set
$$
\Phi_{e_1}:=\left(
\begin{matrix}
1 & 0 & 0 & 0 & 0\cr
0 & 1 & 0 & 0 & 0\cr
0 & 0 & 1 & 0 & 0\cr
1 & 1 & 0 & 1 & 0\cr
0 & 1 & 1 & 0 & 1
\end{matrix}\right),
\quad
\Phi_{e_2}:=\left(
\begin{matrix}
1 & 0 & 0 & 0 & 0\cr
0 & 1 & 0 & 0 & 0\cr
0 & 0 & 1 & 0 & 0\cr
0 & 0 & 1 & 1 & 0\cr
1 & 1 & 0 & 0 & 1
\end{matrix}\right),
\quad
\Phi_{e_3}:=\left(
\begin{matrix}
1 & 0 & 0 & 0 & 0\cr
0 & 1 & 0 & 0 & 0\cr
0 & 0 & 1 & 0 & 0\cr
0 & 1 & 0 & 1 & 0\cr
0 & 0 & 1 & 0 & 1
\end{matrix}\right).
$$
To give the bilinear map $\Theta$ we set
$$\Theta(e_1,e_1):=(1,1,1,1,0)^t,\ \Theta(e_2,e_2):=(1,1,0,0,0)^t,\
\Theta(e_3,e_3):=(1,0,0,0,0)^t,$$
$$\Theta(e_2,e_1):=(1,0,0,1,1)^t,\
\Theta(e_3,e_1):=(0,1,0,0,1)^t,\
\Theta(e_3,e_2):=(0,0,1,1,1)^t$$
with the convention that the $\Theta(e_i,e_j)$ which are not mentioned are
equal to $0$.
From the presentation (\ref{gr2562}) it can be seen that the 
resulting group $\Z_2^5\rfish_{\Phi,\Theta} \Z_2^3$ is isomorphic to 
${\rm G}(256,2)$.

A mixed ramification structure of type
$([4,4,4]_{16}$ 
on 
$\Z_2^5\rfish_{\Phi,\theta} \Z_2^3$ 
is given by 
$$[((0,1,0,1,0)^t,e_3),((0,0,1,0,1)^t,e_2+e_3),((0,0,0,0,0)^t,e_2)].$$
The conventions are the same as in the example G(256,1).

\end{document}